\documentclass[12pt,thmsa]{article}
\usepackage{amsfonts}
\usepackage{amssymb}
\usepackage{graphicx}
\usepackage{amssymb,amsmath}
\usepackage[all,2cell,dvips]{xy} \UseAllTwocells \SilentMatrices
\newtheorem{teo}{Theorem}[section]

\newtheorem{lemma}[teo]{Lemma}

\newtheorem{obs2}[teo]{Remark}

\newtheorem{tea}{Theorem}[subsection]

\newtheorem{no2}[teo]{Note}

\newtheorem{no3}[tea]{Note}

\newcommand{\Gal}{{\rm Gal}}
\newcommand{\Frob}{{\rm Frob }}

\newcommand{\PSL}{{\rm PSL}}
\newcommand{\SL}{{\rm SL}}
\newcommand{\Symm}{{\rm Symm}}

\newcommand{\PGL}{{\rm PGL}}

\newcommand{\GL}{{\rm GL}}

\newcommand{\F}{{\mathbb{F}}}
\newcommand{\Z}{{\mathbb{Z}}}
\newcommand{\Q}{{\mathbb{Q}}}

\def\timehm{\count31=\time \count32=\count31 \divide\count31 by 60
\number\count31 \multiply\count31 by 60 \advance\count32 by
-\count31 :\ifnum\count32<10 0\fi \number\count32}

\newcommand{\qed}{\hfill\rule{2mm}{2mm}}

\def\ideal#1{<\kern-2pt #1\kern-2pt >}

\begin{document}
\title{{\bf  Automorphy of $\Symm^5 (\GL(2))$ and base change}
}
\author{Luis Dieulefait\thanks{Research partially supported
by MICINN grants MTM2009-07024 and MTM2012-33830 and by an ICREA
Academia Research Prize}
\\
Universitat de Barcelona\\
e-mail: ldieulefait@ub.edu\\
%Date:  (Preliminary Version)
 }
\date{\empty}

\maketitle

\vskip -20mm

\begin{quote}
\small \center \it a mi hijita Luana:  Bienvenida al tercer planeta!
Ojal\'a que sus maravillas,
 naturales o imaginadas, nunca dejen de sorprenderte.
\end{quote}

\begin{abstract} We prove that for any Hecke eigenform $f$ of level $1$ (i.e., an eigenform on the complex upper-half plane with respect to $\SL_2(\mathbb{Z})$) and arbitrary weight there is a
self-dual cuspidal
 automorphic form $\pi$
 of $\GL_6(\mathbb{A}_\Q)$ corresponding to $\Symm^5 (f)$, i.e., such that the system
 of Galois representations attached to
 $\pi$ agrees with the $5$-th symmetric power of the one attached to
 $f$. \\
 %Assuming a slight strengthening of the Automorphic Lifting
 %Theorems that we use in the proof of this result, we note that the
 %same conclusion applies also to any newform without CM of level prime to $30$.
 %\\
We also improve the base change result that we obtained in a
previous work: for any newform $f$, and any totally
real number field $F$ (no extra assumptions on $f$ or $F$), we prove the
existence of base change
relative to the extension $F/\Q$. \\
Finally, we combine the previous results to deduce that base change
also holds for $\Symm^5(f)$: for any Hecke eigenform $f$ of level
$1$ and any totally real number field $F$, the automorphic form
corresponding to $\Symm^5 (f)$ can be base changed to $F$.
% (again,
%this can be {\it conditionally} extended to non-CM newforms of level
%prime to $30$).
\end{abstract}
%MSC: 11F80, 11F11, 11R39

\section{Introduction}

After completion of our previous work on base change for $\GL(2)$
(cf. [Di12a]), while explaining the result at a few conferences
we stressed that ``by modifying a bit the proof it should follow
that every pair of newforms can be connected to each other by a {\it
safe chain}". Moreover, we predicted that such a result should have
important consequences, since it should allow to propagate {\it nice
modularity properties} (such as base change) from a single newform
to the rest of them.\\
In this paper, we give a detailed proof of the existence of such
safe chains and we give the first applications of this to Langlands
functoriality. \\
The definition of {\it safe chain} is context-dependent, in all
cases it refers to a series of congruences between two different
newforms (it is actually better to work with the equivalent notion
of congruences between the attached Galois representations), in
({\it a priori}) arbitrary characteristics, having as initial and
final elements the two given newforms. Actually, the congruences are
up to twist by finite order characters, and in this paper (reusing
an idea already appearing in [Di09]) we also allow replacing a
newform by a Galois
conjugate of it as a valid {\it move}.\\
The important thing about safe chains is that (at each congruence in
the chain) they should preserve automorphy of some derived objects:
{\it soit} the restrictions to the Galois group of some number field
of the given pair of modular Galois representations, {\it soit} some
symmetric power $\Symm^n$ of them. By ``preserving automorphy" we
mean that some available Automorphy Lifting Theorem (A.L.T.) should
apply, so as to ensure that IF one of the representations
($2$-dimensional or $n+1$-dimensional, depending on the case) is
automorphic, so is the other. Sometimes (as in [Di12a]) this is
needed to work only in one of the two directions (and this can be
thought as some sort of inductive process propagating modularity
from objects of small ``invariants" to those of higher
``invariants"), but in this paper we will need this to work in both
directions (in order to apply a transitivity argument).\\
Therefore, the conditions imposed at each link (i.e., congruence) in
the chain that makes it {\it safe} depend on the A.L.T. available:
typically we should impose some local condition (specially at the
residual characteristic $p$) at both sides of the congruence, and we
should ensure that the residual image is sufficiently large. \\
We allow ``congruences up to twist" because the effect of twisting a $2$-dimensional
 Galois representation on its $n$-th symmetric power is again a twist, and twists are known to
preserve automorphy, and they are harmless because all the technical
conditions required in the A.L.T. that we will apply (local
conditions, size of residual image) are known to be preserved by
twisting by a finite order character. One of the standard uses of twisting (used for example in [KW]) is to reduce the Serre weight of a residual mod $p$ representation to the range $2 \leq k \leq p+1$ by twisting by a suitable power of the mod $p$ cyclotomic character: we will use this trick repeatedly (observe that such a twist does not affect the Serre's level, i.e., the prime-to-$p$ part of the conductor, of the residual representation).
We also allow Galois
conjugation because the $n$-th symmetric powers of two conjugated Galois representations
are themselves conjugated to each other and Galois conjugation obviously preserves automorphy.\\

 The main result in this paper concerns the application of this
machinery to deduce automorphy of $\Symm^5 (f)$ for all level $1$
cuspforms $f$. The result follows from the construction of a
suitable {\it safe chain}, it suffices to prove two things:\\
1) base case: there exists a level $1$ cuspform $f_0$ such that
$\Symm^5 (f_0)$ is automorphic.\\
2) safe chain: every pair of level $1$ cuspforms $f$ and $f'$ can be
linked by a chain such that the corresponding chain of
$6$-dimensional representations linking $\Symm^5(f)$ with
$\Symm^5(f')$ is a {\it safe chain} in both directions, in such a
way that via suitable A.L.T. we can conclude that $\Symm^5(f)$ is
automorphic if and only if $\Symm^5(f')$ is automorphic.\\

Observe that for the $6$-dimensional Galois representations to be
residually irreducible (this is part of the requirements for the
application of all A.L.T.) one should avoid working in
characteristics $2, 3$ or $5$ at all steps.\\

For part (1), we will consider a cuspform $f_0$ having, in some
large characteristic $p$, residual image projectively isomorphic to
$A_5$, and we will show how to deduce from this residual automorphy
of $\Symm^5 (f_0)$, and then automorphy via a suitable A.L.T. \\

For part (2), we need to construct a series of congruences linking
each given pair of level $1$ cuspforms, making sure that at all
steps the residual characteristic is greater than $5$ and residual
images are ``large" or even ``$6$-extra large" (even after
restriction to $\Q(\zeta_p)$): see the precise definition of these
notions at the end of this introduction. Thanks to the results in
[GHTT] and in Guralnick's Appendix A (see [G]) to this paper, this
is enough to ensure that the $5$-th symmetric powers of these
representations (restricted to $\Q(\zeta_p)$) have residually {\it
adequate} image, a condition needed to apply the two main A.L.T. in
[BLGGT], theorems 4.2.1 and 2.3.1. In fact, we are relying on a
slight strengthening of theorem 4.2.1 of loc. cit. which is proved
in Appendix B (written jointly with T. Gee, see [DG]). Thanks to
this improvement the assumption\footnote{just for this line, we use
the notation of [BLGGT]: the residual characteristic is $\ell$, and
the dimension of the Galois representations is $n$} $\ell \geq
2(n+1)$ (which is part of condition (4) of this theorem, cf.
[BLGGT]) can be replaced by the standard condition of adequacy of
the (restriction of) the residual image, as
in condition (5) of theorem 2.3.1 in loc. cit.\\
We also need to make sure that at each step the $5$-th symmetric
power of both $p$-adic Galois representations being considered
satisfy (locally at $p$) the local condition needed to apply some
A.L.T. in both directions. We will see that in our process they are
either both ``potentially diagonalizable" or both ordinary, allowing
us to apply theorems 4.2.1 (its variant in Appendix B) and 2.3.1 in loc. cit., respectively. \\
The reader should also notice that these A.L.T. are stated in loc.
cit. and in Appendix B for representations of imaginary CM fields,
but a standard argument using quadratic base change allows an
extension to the case of totally real fields, where representations
are now essentially self-dual and automorphy is stated in terms of
 RAESDC automorphic representations, see for example  [CHT]. \\
The fact that our chains are reversible, i.e., that they can be used
to propagate automorphy in both directions, is crucial for us, since
the way that we proceed to link any pair of level $1$ cuspforms is
through an inductive process, combined with transitivity. More
precisely, we are able to link any pair of level $1$ cuspforms $f$
and $f'$ by linking both of them to two newforms in the same orbit,
i.e.,  to two conjugated newform $g$ and $g^\sigma$. Since Galois
conjugation is a valid move (because it preserves automorphy), this
implies that the two given newforms are linked to the same $g$, thus
by transitivity they are linked to each other: just concatenate the
two chains! Thus, if automorphy is known for $\Symm^5 (f')$, it is
propagated to $\Symm^5(f)$ passing through $\Symm^5(g)$ somewhere in
the chain. It is clear then that A.L.T. should work in both
directions, since we need them to propagate automorphy (starting at
$f'$) ``down" until we reach $g$, and
then ``up" until we reach $f$.\\
Let us stress that the orbit of $g$ is unique and independent of the
given forms $f$ and $f'$: it is a universal step appearing in all
safe chains. The inductive argument linking any level $1$ cuspform
$f$ with $g$ will be described in detail in the paper, and it
consists of a series of safe congruences that finally lead to a
space of newforms of small level and weight, and fixed (up to
conjugation) nebentypus, where
computations show that there is a unique orbit of newforms.\\
The construction of this safe chain borrows some key ideas from
[Di12a] and [Di09], but they have to be modified and generalized in
order to be applied in this new context.\\

Using this safe chain, together with the techniques of
``ramification swapping" (as in [Di12a], also used in [Di12b]) and ``killing
ramification" (as in [Di12a]), we will also be able to obtain an
improved proof of base change for $\GL(2)$. This time the safe chain
will connect any cuspform to a CM form, and now ``safe"
will refer to the fact that Modularity Lifting Theorems (M.L.T.) for
the restrictions of the $2$-dimensional Galois representations to
the Galois group of a totally real number field do apply in both
directions at each step of the chain. \\
We are relying again, for the ramification swapping process,  on the
two main A.L.T. in [BLGGT] and the improvement in Appendix B, this
time for $2$-dimensional Galois representations, and the important
feature of these theorems that the required local conditions are
preserved under arbitrary base change allows us to work over
arbitrary totally real fields $F$, not needing to impose any local
condition. For the ``killing ramification" process, we will need to
apply the M.L.T. of Kisin in [K-1], which requires the residual
characteristic $p$ to be split in $F$, but since this step will take
place at a point where the prime $p$ to be killed is a suitable
auxiliary prime, we just have to check that such primes can be
chosen to be split in $F$ (a similar process took
place in [Di12a]).\\
In order to cover the case of newforms of even level we will also apply a $2$-adic M.L.T.
of Kisin in [K-2] together with some ideas taken from [KW].\\

We also prove a base change result for the automorphic forms
corresponding to $\Symm^5(f)$, which follows easily by combining
previous results.\\

We conclude this introduction by warning the reader that the proof
of the main result (the construction of the safe chain) is full of
twists and turns, and at several places the path is so full of
obstacles that it seems there is no way to continue, and it is only
by some completely new trick, or by an unexpected
generalization/combination of results in some of our previous works
that one manages to keep going. It is a road full of miracles,
including a fascinating generalization of Maeda's conjecture that
motivates the choice of the ``small single orbit space" of newforms
at the bottom of the road. The tricks of ``good-dihedral primes",
``micro good-dihedral primes" (for short, MGD primes), ``Sophie Germain primes", ``weight
reduction via Galois conjugation", are
basic to build the safe chain.\\
Both appendices are also of key importance: without them we would
have been forced to work only on characteristics $p
>13$, and just from the computational point of view (to say nothing
of the theoretical complications) it would have been impossible to
complete
the ``low part" of the chain. \\

Let us now record the two main theorems proved in this paper: in sections 2,3, and 4,
 and in section 6 for the extension to other fields, we prove:

\begin{teo} Let $f$ be a level $1$ cuspform. Then if we call $\rho$ an $\ell$-adic Galois
representation attached to $f$, the $6$-dimensional Galois representation
 $\Symm^5(\rho)$ is automorphic, namely, there is a cuspidal self-dual automorphic
 representation $\pi$ of $\GL_6(\mathbb{A}_\Q)$ (a RAESDC in the notation of [CHT])
such that the  $\ell$-adic Galois representation attached to $\pi$ is isomorphic
 to $\Symm^5(\rho)$.\\
Moreover, if we restrict this $6$-dimensional Galois representation to any totally
real number field $F$, it is still automorphic, meaning that it is isomorphic to
the $\ell$-adic Galois representation attached to a cuspidal self-dual automorphic
 representation of $\GL_6(\mathbb{A}_F)$.

\end{teo}

%Furthermore, a conditional extension of this result to non-CM newforms of  level
 %prime to $30$ is discussed in section 7. \\

In section 5 we prove base change for classical cuspidal modular forms in full generality:

\begin{teo} Let $f$ be a newform of  arbitrary level and weight $k \geq 2$. Let $F$ be any
 totally real number field. Then, $f$ can be lifted to $F$, i.e., there is a
  Hilbert newform $f'$ over $F$ whose attached $\ell$-adic Galois representations
   are isomorphic to the restriction to $G_F$ of those attached to $f$.
\end{teo}

{\bf Remark}: The (safe) chain constructed in this paper connecting
any pair of newforms to each other has been applied in a sequel to
this paper to prove automorphy of $m$-fold tensor products of level
$1$ cuspforms for any $m$ (assuming regularity of the tensor
product).
See [Di12c] for details. \\

{\bf Acknowledgments:} My greatest thanks are due to R. Guralnick
for his immediate response to my S.O.S. asking for better results on
adequate groups, and for writing these results in Appendix A,  and
to R. Taylor for kindly helping me through discussions to check some
parts of this proof, and for kindly having me as guest at the I.A.S.
for a period of three weeks during which some of the most
important advances in this proof where made. \\
I also want to thank T. Gee for providing the key insights needed
for the proof of the result in Appendix B (and for pointing out a
slight inaccuracy in a previous version of this work),
 P. Tsaknias for some vital computations that
 he was kind enough to do for me on MAGMA, and the Number Theory Group
  of The Department of Mathematics of
the University of Luxembourg for providing him the required
computers. For some useful conversations, remarks, explanations, and
corrections on previous versions, I also want to thank
 M. Dimitrov, J.-P. Serre, M. Kisin, Bao Le Hung and F. Herzig. \\
 Finally, I want to give special thanks to the referee, whose
 suggestions have greatly improved the presentation and readability
 of this paper.\\

We finish this section with some definitions and notations.\\

\bf{Notation}: \rm In this paper, $F$ will always denote a totally real number field. \\
For every number field $K$, we will denote by $G_K$ the absolute Galois group of $K$. \\
We will write $\zeta_p$ for a primitive $p$-th root of unity.
 We
will denote by $\chi$ the $p$-adic or mod $p$ cyclotomic character.
 The value of $p$, and whether it is the $p$-adic or the mod $p$ character, will
 always be clear from the context.\\
We will denote by $\omega$ a Teichmuller lift of the mod $p$ cyclotomic character.\\
Given a Galois representation $\sigma$, we will denote by $\mathbb{P}(\sigma)$ its projectivization.\\

\bf{Definitions}: \rm Let $K$ be a number field. Let $\bar{\rho}_p$
be a two-dimensional, odd, representation of $G_K$ with values on a
finite extension of $\F_p$.
\begin{enumerate}  \item We say that the image of
$\bar{\rho}_p$ is \it{large} \rm if $p \geq 7$ and the
 image contains $\SL(2, \F_p)$, or $p = 3$ or $5$  and the image contains
  $\SL(2, \F_{p^r})$ for some $r \geq 2$.
 If this is the case, it is easy
 to see that the image of $\mathbb{P}( \bar{\rho}_p )$ is
 isomorphic to one of the following two groups: $\PSL(2, \F_{p^r})$
 or $\PGL(2, \F_{p^r})$, for some $r$ (and $r \geq 2$ if $p=3, 5$). This implies
 in particular that large images are non-solvable, and also that
 they are adequate subgroups of $\GL(2, \bar{\F}_p)$, see theorem 1.5 and
  remark 1.6 (3) in Appendix A and the main result of [GHTT].

 \item We say that the image of
$\bar{\rho}_p$ is \it{$6$-extra large} \rm if $p = 11$ or $13$ and
the
 image contains $\SL(2, \F_{p^r})$ for some $r \geq 2$. \\
 For $p=7$ or $p \geq 17$ we say that the image of
$\bar{\rho}_p$ is \it{$6$-extra large} \rm if it is large. \\
 This technical condition is
 needed to ensure that the $5$-th symmetric power of $\bar{\rho}_p$
 has adequate image (for $p \geq 17$ this is just an application of the main result in [GHTT], for $p=7$ it follows from remarks 1.6 (1) and 1.6 (3) in Appendix A, and for $p = 11$ and $p=13$ it follows from theorem 1.5 combined with remark 1.6 (3) in Appendix A). \\

{\bf Remark}: Observe that since large and $6$-extra large images correspond to
almost simple projective images,  these conditions are preserved if we restrict
the representation to a
cyclic extension of $K$, such as $K(\zeta_p)$.\\

    \item We say that the image of $\bar{\rho}_p$ is \it{dihedral} \rm when
     the image of $\mathbb{P}( \bar{\rho}_p )$ is a dihedral group of order at least $4$.

    \item We say that the image of $\bar{\rho}_p$ is \it{bad-dihedral} \rm when
    it is dihedral, $p >2$, and  the restriction
    of $\bar{\rho}_p$ to the absolute Galois group of $K(\sqrt{\pm p})$ becomes reducible,
     where the sign is $(-1)^{(p-1)/2}$.

\end{enumerate}

\section{Construction of the safe chain: a general description}
This is the {\it code description} of what is needed to make the
safe chain. Recall that this is just ``half" of the chain, the full
chain is obtained just by concatenation of two of these safe chains,
which makes sense because both end up in the same orbit, and are
reversible.\\
At each step we describe the nature of the movements to be
performed, and the output. References are to indicate where some of these tricks first appeared.\\
 We start with a cuspform $f$ of level $1$ and weight $k \geq 12$: \\

1) Introduce Good-Dihedral prime [KW]: level raises to $q^2$. \\
2) Weight reduction via Galois conjugation [Di09]: weight reduces to $2 < k \leq 14$ ($k$ even).\\
 3) {\it Ad hoc} tricks to
make the small weight congruent to $2$ mod $3$ (Sophie Germain
primes, Hida families...): end up with $k \equiv 2 \pmod{3}$ and $k
< 43$.
\\
4) Introduce MGD prime $43$ using the pivot primes $7$ and $11$
[Di12a]: end up with newform of weight $2$ and
level $43^2 \cdot q^2$.\\
5) Remove the Good-Dihedral prime (in two moves): end up with a newform of weight $q+1$ and level $43^2$.\\
6)  Again weight reduction via Galois conjugation but this time
``highly improved", because we need to ensure large residual image
at each step using just the MGD prime.
Weight reduces to  $2 < k \leq 14$ ($k$ even).\\
7) {\it Ad hoc} tricks to make weight smaller than $17$ and divisible by $4$
 (Sophie Germain primes, Khare's weight
reduction, non-dihedrality due to class
field theory of real quadratic field...): end up with a newform of weight
 $16$ and level $43^2$. \\
8)  Introduce nebentypus at $17$ of order $8$: end up with a newform of
 weight $2$, level $43^2 \cdot 17$, with nebentypus. \\
9)  Remove the MGD prime modulo $11$ via an {\it ad hoc} Lemma to
ensure residual irreducibility: get congruence (maybe using
level-raising)
 with a newform Steinberg at $43$, of weight $2$ and level $43 \cdot 17$,
 with nebentypus of order $8$ at $17$. \\
10) Move from weight $2$ to weight $44$ by reducing modulo $43$: irreducibility checked by hand. \\
11)  Generalized Maeda in weight $44$, level $17$, nebentypus of
order
$8$: check that this space has a unique orbit.\\

Several steps use tricks appearing in the referred paper, and they
will be easy to follow for those who have read these papers. Steps
(3) and (7) are painful and technical: the weight has been made
small, BUT we require some technical condition on it to be able to
perform the next step (i.e., steps (4) and (8)): step (3) is done
because at step (4) we want $k \equiv 2 \pmod{3}$ because the prime
$43$ is not truly a Sophie Germain prime, there is a $3$ dividing $43-1$,
thus when we take a modular weight $2$ lift in characteristic $43$
IF WE KNOW that $k \equiv 2 \pmod{3}$, the nebentypus, which is
$\omega^{k-2}$, will be a character of order $7$, and this allows us
to work with $43$ ``as if it were a Sophie Germain prime", and using
this character of order $7$, and later the factor $11 \mid 43+1$, we
manage to introduce $43$ as a MGD prime (compare with the
introduction of $7$ as a MGD prime in [Di12a]). Step (7) is done
because before moving to characteristic $17$ in step (8) with a
weight $k$ smaller than $17$ to take a modular weight $2$ lift,
since we want the character $\omega^{k-2}$ to be of order $8$, we
NEED to ENSURE that $k$ is divisible by $4$. \\
Both in steps (3) and (7), in order to ``reduce" to weights
satisfying the required conditions, we make a few technical moves
that allow ``magically" to control the weight and force it to
satisfy them. \\
Step (6) is perhaps the best technical innovation in the paper: we
need to perform the weight reduction to inductively go from an
arbitrarily large weight $k$ to $k \leq 14$, we do this using the
method of [Di09] as in step (2), but this time the level contains
just the MGD prime 43, and still we need to ensure that in the whole
process of weight reduction all residual images are going to be
large! We succeed thanks to an extra degree of freedom, previously
unexploited, in Galois conjugations: there is not a unique
conjugation that reduces the weight, there are at least two
completely different choices (this we prove), and we can show that
for at least one of them the residual representation will have large
image. \\
Step (9), if it were performed by direct computation, will not be
surprising. But since the level at this step is $17 \cdot 43^2$ and
the nebentypus of order $8$, computations seem out of reach. We prove
irreducibility modulo $11$ (recall that ramification at $43$ being
of order $11$, at this step one is losing the MGD prime from the
level) by using in a non-trivial way information on the ramification
at $17$ and on the trace at $43$ of the residual representation,
assumed reducible, to derive a contradiction.\\
The final step (11) is a direct computation to check that certain
space has a unique orbit, but this does not come out of the blue: we
``knew" that there should be a unique orbit there by a precise
generalization to arbitrary level of Maeda's conjecture (cf. [DT]).\\

\section{The eleven steps of the safe chain}

Before describing each step of the chain, let us make some general
remarks. \\
First of all, let us stress that in all steps we are always working
in residual characteristic $p \geq 7$. This is necessary because we
want the symmetric $5$-th power to be residually
irreducible (any $2$-dimensional representation over a field of characteristic up to $p$ has a reducible symmetric $p$-th power).\\
 At each step, we will introduce several congruences
between modular Galois representations. Since we want to apply the
two main A.L.T. in [BLGGT], Theorems 4.2.1 (its variant in Appendix
B) and 2.3.1, in both directions, to the $5$-th symmetric powers of
these representations, at each of the congruences introduced we need
to ensure that both $6$-dimensional representations are
``potentially diagonalizable" (PO-DI) or that they are both
ordinary. Observe that both of these local properties are preserved
by taking symmetric powers, so we just have to check them at the
level of the modular $2$-dimensional Galois representations (for the
PO-DI property, just observe that if two points in a local
deformation ring $R_{\mathcal{D},p}$ are in the same irreducible
component, when we switch to a larger deformation ring
$R_{\mathcal{D'},p}$, i.e., with the latter ring projecting into the
former one, the corresponding points in the larger ring are also in
the same irreducible component: this property is already applied
in [BLGGT], page 531, to show that this notion is preserved by base change, i.e., when one restricts
the Galois representations. Of course, we are also using the elementary fact that a symmetric power of a sum of crystalline characters is again a sum of crystalline characters). This will also be useful when we
address, in latter sections, the problem of base change, since for
that case we are going to rely on the same A.L.T. most of the time,
but this time applied to $2$-dimensional Galois representations.\\
Most of the congruences that we will perform involve taking weight
$2$ modular lifts, for a residual representation in some odd
characteristic $p$ having Serre's weight $2 < k < p+1$, and the
congruence will be between a potentially Barsotti-Tate
representation and a crystalline representation in the
Fontaine-Laffaille range (i.e., of weight $k \leq p$). Since in both
cases (by the results of Kisin in [K-BT] in the first case, as
proved in [GK], and by [BLGGT] and [GL] in the other) such
representations are known to be PO-DI, the local conditions
needed to apply Theorem 4.2.1 in [BLGGT] and its strengthening in Appendix B are satisfied. \\
The second (and last) case that we will encounter in our chain is
when one of the representations is semistable at $p$ of weight $2$,
and the other one is crystalline at $p$ of weight $p+1$. In this
case, both representations are known to be ordinary (in fact, they
live in the same Hida family), thus again we can apply  an A.L.T. in
[BLGGT], this time Theorem 2.3.1, the one for ordinary
representations (see [G] for the definition of ordinary
representations). At one step (step 3) there is a small variant of
this: we consider a congruence between a weight $2$ semistable
representation
 and a crystalline representation whose weight is actually larger than $p+1$:
 we take a suitable specialization of
  the Hida family at a larger weight $k \equiv 2 \pmod{p-1}$. Again, this is
   know to be an ordinary representation,
   so this A.L.T.  applies.\\
The reader can check (this is automatic) that in the $11$ steps that
follow we will always be in one of the two cases above, thus from
now on we are not going to insist anymore on this point: the local
conditions required to apply the A.L.T. to the $5$-th symmetric
powers of our representations are satisfied. Also, when in latter
sections we discuss base change, since these local properties are
preserved by base change, both for the $2$-dimensional
representations and for their symmetric powers, the local conditions
to apply the A.L.T.
over any totally real number field $F$ are satisfied.\\
The other obstacle to apply the A.L.T. in [BLGGT] (recall once again
that concerning theorem 4.2.1 in loc. cit., we are relying on the
improvement that we prove on Appendix B, cf. [DG]) is the condition
that residual images should always be adequate. We will check that
in our chain residual images will always be $6$-extra large (this is
not so at the ``base case", but in this case the image is known to
be adequate, too), a condition that is preserved after restriction
to the Galois group of $\Q(\zeta_p)$, and this implies that the
images of these restrictions are adequate, thanks in
particular to the results in Appendix A. \\
The required condition on the residual images will be obtained, at
most steps, by using some ramification information. This will be
particularly easy when a Good-Dihedral prime is in the level, and at
the very bottom of the chain, when no significant ramification is
preserved, will be checked partly by hand.\\
In any case, let us stress that in what follows we are just required
to
indicate two things at each step:\\
i) how the chain is built, and \\
ii) how it is checked that the residual images are $6$-extra large
\\
 Observe that since
 $6$-extra large implies large, when we re-use this chain for base
 chain, we at least know that it is built with large residual images
 (and we will just have to check that this largeness is preserved by
 base change).\\

\subsection{Step 1}
We start with $f$ of level $1$ and weight $k \geq 12$. At this step
we introduce a Good-Dihedral prime. This is done just as in [Di12a]
(and also in [KW] where the idea first appeared), so we will be
brief. We first choose a prime $r > 13$ larger than $k$ where the
residual image is large, thus $6$-extra large, and reducing modulo
$r$ we switch to a weight $2$ situation (this can be done for
example by applying theorem 5.1- (2) in [KW] combined with a
suitable M.L.T., such as theorem 4.1 in loc.cit.). Here we are
introducing a nebentypus at $r$, to be removed later. More
precisely, we have now a newform in $S_2(\Gamma_1(r))$ with
nontrivial nebentypus (this more refined information also follows
from theorem 5.1 - (2) in [KW]) . Then, we take primes $t$ and $q$
larger than $r$ as in Lemma 3.3 of loc. cit., and as in the
discussion thereafter: by working modulo $t$, where the image is
$\GL(2, \F_t)$, thus $6$-extra large, we produce a congruence with a
newform $f_2$ of weight $2$ with $q^2$ in the level, such that
ramification at $q$ is given by a character of order a power of $t$,
and $t$ divides $q+1$. Let us stress that there exist infinitely
many primes $q$ satisfying the conditions of Lemma 3.3 of loc. cit.
The prime $q$ is a supercuspidal prime for $f_2$. Before choosing
$t$ and $q$, a large bound $B$ is fixed, greater than $k$, $ 2 \cdot
r$, and than any other auxiliary prime that one will use (such as
$p= 11, 43$) and $t$ and $q$ are chosen to be larger
than $B$. We take $B  > 68$ because we are going to need this in Step 3.\\
These two primes are also asked to satisfy:\\
$t \equiv 1 \pmod{4}$, $q \equiv 1 \pmod{8}$, and also, for every
prime $p \leq B$, we require that $q \equiv 1 \pmod{p}$.\\
In the following steps, up to step 4 (included), we are going to
work all the time in characteristic $p \leq B$, and this implies
that the local ramification at $q$ will be preserved through these
steps. This, together with the above conditions on $q$ implies, and
this is the main point about Good-Dihedral primes (see ``A very important remark" in page 1023 of [Di12a]), that the residual images will be large. Moreover, since we
are taking $t > B > 2 \cdot r > 13 $ and the residual projective
images contain an element of order $t$ (given by the image of the
inertia group at $q$) it is obvious that if we are in characteristic
$11$ or $13$ (and trivially in any other characteristic $7 \leq p
\leq B$)
the residual image will also be $6$-extra large. \\
We conclude the first step by moving back to characteristic $r$. We
take a mod $r$ Galois representation attached to $f_2$ and we lift
it to a modular representation attached to a newform $f_3$ without
$r$ in the level. As explained in the introduction, we freely twist
when necessary, so we are assuming that the Serre's weight of this
residual representation is at most $r-1$ (and greater than $2$
because the nebentypus of $f_2$ at $r$ was not trivial), and $f_3$
is taken to be of this weight. Therefore, we end up with $f_3$ of
level $q^2$, Good-Dihedral at $q$, and weight $2< k < B/2$.\\

At a referee's suggestion, we produce here a picture of the congruences just described linking the given $f$ with $f_3$. The superscript $q$-GD means that the prime $q$ is a Good-Dihedral prime for a newform, and the symbol $=^*$ denotes that two residual representations are isomorphic after twisting by some power of $\chi$:

$$
\xymatrix{ (\rho_{f_1, p} )_p  \ni     \rho_{f_1, t} \ar[rd]  &     &  (\rho_{f_2, p} )_p  \ni     \rho_{f_2, r} \ar[rd]    &   (\rho_{f_3, p} )_p \\
   f  \equiv_r  f_1  \in  S_2(\Gamma_1(r)) \ar[u] & \bar{\rho}_{f_1, t}   \ar[ru]^{KW \atop  Thm \; 5.1-(4)}  &   f_2   \in  S_2^{q-GD}(\Gamma_1(q^2 r)) &    \bar{\rho}_{f_2, r}  =^* \bar{\rho}_{f_3, r} , f_3 \in S_{k<r}^{q-GD}(\Gamma_1(q^2)) \ar[u]}
$$

\subsection{Step 2}
At this step we perform the Weight Reduction via Galois Conjugation
(WRGC), exactly as in [Di09] (except that when the weight is $10$,
we are not reducing it because we don't need to), so as to reduce to
the case $k \leq 14$.\\
As explained in the previous step, since we are going to work in
characteristics $p  \leq B$ and we have a large Good-Dihedral prime in
the level, residual images are all going to be $6$-extra large.\\
Observe that since $k < B/2$, we are allowed to work with primes
larger than $k$ as long as they satisfy $p < 2 \cdot k < B$, and
this is the case in the process of WRGC.\\
We also note for the reader's convenience that, together with Galois
Conjugations (an allowed move, see the introduction), this weight
reduction involves congruences between potentially Barsotti-Tate and
Fontaine-Laffaille type Galois representations, exclusively.\\
Since the process is described in [Di09], we are going to be very
brief:\\
Suppose $k >14$ (if $2 < k \leq 14$ we can go directly to step 3).
It is perhaps worth recalling that through all this process weights
are
going to be even, and always larger than $2$.\\
Let $p$ be the smallest prime larger than $k$, except in the case $k
= 32$
where we take $ p = 43$.\\
Consider the mod $p$ Galois representation attached to $f_3$, and
take a modular weight $2$ lift of it, corresponding to a newform
$f_4$ having nebentypus at $p$ given by $\omega^{k-2}$, a character
of order $m = (p-1)/d$, where $d= (p-1, k-2)$. As in Step 1, this is done by applying theorem 5.1 - (2) of [KW] (in combination with a suitable M.L.T., such as theorem 4.1 in loc. cit.).
Let us call $a = (k-2)/d$.\\
The level at $p$ of $f_4$ is $p$ and it is principal series at $p$,
with inertia Weil-Deligne parameter of the attached $\ell$-adic Galois representations
at $p$ (and the same for the $p$-adic Galois
representation) being given by $\omega^{k-2} \oplus 1$.\\
 Let us assume (this will be proved later)
that $m
> 6$, and let us consider $t$ to be the smallest integer greater than $m/2$
and relatively prime to $m$. \\
Warning: this integer has no relation with the prime $t$ giving the order of ramification at the Good-Dihedral prime $q$ used in Step 1. We hope that the use of the same letter to denote both objects will not be a cause of confusion because they live in different worlds.\\
Observe that the
difference $t - m/2$ is at most $2$, and that $t < m-1$. \\
The Galois conjugation is meant to change the nebentypus: we take an
element $\gamma$ in the Galois group of $m$-th root of unity
corresponding to raising to the $a^{-1} \cdot t$, where the inverse
of $a$ is taken modulo $p-1$. This changes the nebentypus (and the
local type at
$p$) from $\omega^{k-2} = (\omega^d)^a$ to $(\omega^d)^t$. \\
Now we consider the mod $p$ Galois representation attached to the
conjugated newform $f_4^\gamma$. Its Serre's weight can be shown to
be (after suitable twisting) either $k_1 = dt + 2$, or $ k_2= p + 1 - d t$, two values whose
sum is $p+3$  (cf. [Sa], Corollary 6.15).
\\
Let us observe that $k_2 < k_1$: this inequality, appearing already in [Di09], is due
 to the fact that $k_1 > (p+3)/2$, which follows easily from the fact that we are taking $t > m/2$. \\
The idea is that because of a strong version of Bertrand's postulate
$p$ will be ``near" $k$, and this will imply that $k_1$ and $k_2$
are smaller than $k$: in fact, since $k_2 < k_1$,
 it is enough to check this for $k_1$, and since $k_1$ is
approximately $d \cdot t$, and this is approximately $(d \cdot
m)/2$, and this is approximately $p/2$, it is clear that having $p$
near $k$
the new weight $k_1$ or $k_2$ will be smaller than the given $k$.\\
More precisely: for $k \geq 38$, it is known that if we consider the
primes nearest $k$: $p_n < k < p_{n+1}$, it holds:
$$ p_{n+1} / p_n  < 1.144 \qquad \quad (*)$$
And also its obvious consequence:
$$ (p_{n+1} - 1) / (p_n -1)  < 1.15 $$
From this, calling again $p = p_{n+1}$ the first prime after $k$, we
have:
$$ (p - 1) / (k-2) < 6/5$$
and this implies that the order $m$ of the nebentypus $\omega^{k-2}$
is greater than $6$ (this we had assumed before, now we are proving it for $k > 36$). \\
Having this, we define $t$ as above and Galois conjugate using it as
above and easily deduce that for the any of the two possible Serre's
weights $k' = k_1$ or $k_2$ of the Galois conjugated residual
representation it holds: $p/ k' > 1.144$.\\
In fact, since $k_2 < k_1$ it is enough to check this for $k' =
k_1$, and then it reads: $p / k' = p / (dt + 2) = (md + 1) / (dt +
2) $. That this is larger than $1.144$ can be easily checked once
the precise value of $t$ is introduced: there are three cases,
depending on whether $m$ is $0$ or $2$ modulo $4$, or odd, but in
any case the above quotients give increasing functions (increasing
with $m$) that already satisfy the inequality for the smallest
values of $m$ ($8$, $10$ and $7$, respectively) and for every $d$,
thus we conclude that the inequality holds in general.\\
Because of $(*)$, this implies: $k' < p_n$, thus in
particular $k' < k$. Thus the induction works.\\
If $14 < k \leq 36$, it can be checked by hand that the above
process also works, exactly as described, and always gives a smaller
Serre's weight $k'$ (and it is easy to see that it always gives even
values, and that they are always greater than $2$).\\
The only subtle point is that for $k=32$ one should take $p=43$.\\
Let us include here, for the reader's convenience, four examples,
including the one of weight $32$:\\

$\bullet$ $k = 16, p = 17$: $d = 2, m = 8, t = 5, dt = 10$; thus: $k' = 12$ or $8$.\\

$\bullet$ $k = 18, p = 19$: $d = 2, m = 9, t = 5, dt = 10$; thus: $k' = 12$ or $10$.\\

$\bullet$ $k = 30, p = 31$: $d = 2, m = 15, t = 8, dt = 16$; thus $k' = 18$ or $16$.\\

$\bullet$  $k = 32, p = 43$: $d = 6, m = 7, t = 4, dt = 24$; thus $k' = 26$ or $20$.\\

This concludes the weight reduction, by iterating the above procedure we end up this step with a
newform $f_5$ of level $q^2$, Good-Dihedral at $q$, of weight $2 < k
\leq
14$.\\
Remark: the inequality $k>2$ follows from the fact that the nebentypus at $p$ of the newforms we are considering is non-trivial. In fact, observe that from the formula for $k_1$ and $k_2$ that we are using it follows that if any of them were equal to $2$ then the nebentypus must be $\omega^0$ or $\omega^{p-1}$, i.e., the trivial character.\\

At a referee's suggestion, we include a picture describing the congruences involved in this step (see end of section 3.1 for notation):

$$
\xymatrix{  f_3 \in S_{k>14}^{q-GD}(\Gamma_1(q^2)) \ar[d]_{\mbox{suitable} \atop p>k} & (\rho_{f_4, p} )_p  \ni     \rho_{f_4, p}  \ar[r]^{\mbox{conjugate}} &  \rho_{f_4, p}^\gamma  \ar[d]     \\
  \bar{\rho}_{f_3, p} \ar[ru]^{KW \atop Thm \; 5.1-(2)}  &   f_4 \in    S_2^{q-GD}(\Gamma_1(p q^2))    &    \overline{ \rho_{f_4, p}^\gamma }  =^*  \bar{\rho}_{f_?, p},  \;  f_?  \in  S_{2 < k' < k}^{q-GD}(\Gamma_1(q^2))
}
$$

If the weight $k'$ of $f_?$ is larger than $14$, we repeat the process until we end up with $f_5$ having weight  $2 < k \leq 14$.

\subsection{Step 3}
As it will be clear in Step 4, before introducing the
Micro-Good-Dihedral (MGD) prime $43$ in the level, we need to reduce
to a situation in which the weight not only is small (smaller than $43$ will be enough)
 but also is congruent to $2$ modulo $3$. Here we combine
some ad hoc tricks that allow us to do so.\\
The ideas in this step are the following: (in logical order, which,
as usual, does not coincide with the final ordering on the paper): \\
If one happens to have some weight $k$ satisfying $k \equiv $2$
\pmod{3}$ and you believe that there is a prime $p >k$ ``near" $k$
such that $p \equiv 1 \pmod{3}$, then you can check (it works in the
few applications that we are going to do now, because it always
works!) that if you do weight reduction using $p$, as in Step 2, the
new (smaller) values of the Serre's weight, $k_1$ and $k_2$, will
both be again congruent to $2$ modulo $3$ (it makes sense because
they sum up $p+3$ and this is congruent to $1 \pmod{3}$).\\
Thus, the strategy is the following: we use the Sophie Germain
trick, taken from [Di12b], this allows you to, given a pair $p_1,
p_2 = 2 p_1 + 1$ of Sophie Germain primes, reduce all cases of even
weights $ 2< k < p_2 $ to a single case, namely, to the case of weight $p_2 + 1$.
Unfortunately, because we are forced to take $p_1  \neq 3$ (we are
not allowed to work on characteristics $\leq 5$), this is not
satisfactory because $p_2 +1$ will never be congruent to $2$ modulo
$3$, but moving to a weight higher than $p_2 + 1$ (using Hida
families) we will remedy this problem. After this step, we will end
up with a unique weight, it will be congruent to $2$ modulo $3$, and
we will complete the step by reducing it, as indicated above, in
such a way that the mod $3$ property on the weight is preserved.\\
We are going to work  all the time in characteristics $p \leq 67$
(and, as usual, $p >5$) and with weights $k \leq 68$, thus since $68
< B < t < q$ the Good-Dihedral prime in the level will ensure that
residual images
are $6$-extra large at each step.\\
Remark: the value $k=68$ appears here because for the pair of Sophie Germain primes $(p_1, p_2) = (11,23)$ that we are going to use the smallest value of a weight $k >2$ such that $k$ is congruent to $2$ modulo $3$ and congruent to $2$ modulo $p_2 -1$ is $68$ (see below for an explanation of why the second of these congruential conditions is important).\\
 The Sophie Germain pair that we consider is given by $p_1 = 11$,
$p_2 = 23$. As in [Di12b], since the weight of $f_5$ satisfies $2 <
k \leq 14$, we move to characteristic $23$, reduce mod $23$, and
take a modular weight $2$ lift, having a nebentypus $\omega^{k-2}$ ramified at $p_2 = 23$
of order $11$. Then we move to characteristic $11$, reduce mod $11$,
and take another modular weight $2$ lift, where we have changed the
type at $23$. Since ramification at $23$ was given by a character of
order $11$, we know that the residual representation will either be
unramified at $23$ or will have semistable (i.e., unipotent)
ramification at $23$. Moreover, in the first case it can be checked
(a similar situation can be found in [Di12a]) that Steinberg
level-raising at $23$ can be performed. It is important to remark that this follows from the fact that this residual representation is the reduction of a representation with ramification at $23$ given by a character of order $11$. Thus we conclude that in any
case we can take a modular weight $2$ lift of this mod $11$
representation that it
semistable at $23$ (i.e., the ramification group at $23$ is pro-unipotent), which corresponds to the fact that the modular form is Steinberg at $23$ (this is either an application of Ribet's level raising to add the prime $23$ or a minimal lift). \\
In the usual Sophie Germain trick, at this step, reducing modulo
$23$ one manages to create a congruence with another newform of
weight $24$. But since such a weight is not useful for our purposes,
we take another specialization of the ordinary Hida family containing this
semistable weight $2$ representation: we take $k = 68$. Since $68
\equiv 2 \pmod{22}$ we know that there is an ordinary Galois
representation attached to a newform of this weight and level $q^2$
which is congruent modulo $23$ with the given weight $2$
representation.\\
We include a proof of the fact that at weight $68$ (and with the same argument, at all weights $k>2$ such that $k$ is congruent to $2$ modulo $22$) the specialization of this Hida family gives a newform without $23$ in its level. This is well-known, but the referee asked us to include a proof, for example, the one that he provided: let us call $f$ the weight $2$ newform that is Steinberg at $23$ and let us consider the ordinary Hida family containing the $23$-adic Galois representation attached to $f$. We specialize this Hida family to weight $68 = 2 + 3 (23-1)^{2-1} $ to get a form $g$ of tame level $q^2$ (not necessarily exact level $q^2$) such that $\pi_{g, 23}$ is a (possibly ramified) principal series $\pi(\chi_1 , \chi_2)$ where one of the $\chi_i$ has to be unramified or Steinberg. Observe that the central character of $g$ at $23$ is given by $\omega_{23} \cdot \omega_{23}^{1-68} = \omega_{23}^{-66} = 1$. Hence $\pi_{g, 23}$ must be an unramified principal series otherwise the central character would be non-trivial. This proves that $g$ has exact level $q^2$.\\
Remark: we insist on the fact that our main goal is to manage to reduce the set of possible weights to values that are all greater than $2$ and congruent to $2$ modulo $3$. This will be essential on Step 4 because thanks to this property on the weights and the fact that the nebentypus at $43$ will be $\omega^{k-2}$ we see that the nebentypus has prime order $42/ 6 = 7$, and this will be required. For this reason, in this Hida family we skip the value $k= 46$ and we choose $k=68$ instead.\\
At this point, we have reduced the case of weight up to $14$
(greater than $2$) to the case of weight $68$. Since $68 \equiv 2
\pmod{3}$, let us try our weight reduction. \\
We must work only with primes that are congruent to $1$ modulo $3$.
Our first try is $p= 73$, using it we manage to reduce the weight but we get
stuck in the next step (with $k' = 44$: no primes congruent to $1$
mod $3$ show up ``nearby"!), so let us try with $p=79$
instead.\\
We reduce modulo $79$, and we take a modular weight $2$ lift, with
nebentypus
$\omega^{66} = (\omega^6)^{11}$, of order $m= 13$. \\
Galois conjugation is done as in Step 2, using $t= (m+1)/2 = 7$, i.e.,
exponentiating to the $11^{-1} \cdot 7$ in the cyclotomic field of
$13$-th roots of unity, but this leads us to a value of $k_1 = d t +
2 =
6 \cdot 7 + 2 = 44$ where, again, we get stuck.\\
Luckily, this is not the only possible Galois conjugation (this
extra degree of freedom will be exploited more systematically at
Step 6), so let us try with other value of $t$, always prime to $13$
of course. Taking $t=8$ the values we obtain are $k_1 = 50$ and $k_2
= 32$.
Observe that both values are again congruent to $2$ modulo $3$ (as we explained, it
 is easy to see a priori that this will hold, due to our careful choice of the prime $p$).\\
Since the first of these two values is too big (our goal includes
the condition $k < 43$), assuming we end up with $k = 50$ we iterate
the process. We take $p=61$, reduce mod $61$, and this time the
nebentypus of the modular weight $2$ lift is
$\omega^{48} = (\omega^{12})^4$, whose order is $5$. \\
Remark: in Step 2 we always had $m >6$, but $m=5$ is also fine for the Galois conjugation trick
because it also satisfies $\phi(m) > 2$ (compare with  [Di09], where
in the case of $k=10$
 weight reduction was done with nebentypus of order $5$. On the other hand, it is
 easy to see that if $m = 6$ or $m<5$ since $\phi(m) \leq 2$ the method can not work).
\\
 Taking $t=3$ for the
Galois conjugation, we easily compute $k_1 = 12 \cdot 3 + 2= 38$ and
$k_2 = 26$.\\
The conclusion is that we have linked the given newform $f_5$ of
small weight, to another newform with the same level $q^2$,
Good-Dihedral at $q$, whose weight is equal to $26, 32$ or $38$.\\
The three values are congruent to $2$ modulo $3$, and are smaller
than $43$.\\
Remark: the referee asks if this process can be continued to produce weights that are even smaller. The answer is that this is very likely to be the case (according to the referee's computations, we can for example reduce to $k=20, 26 ,32$) but since we already have a set of weights that are all smaller than $43$ and congruent to $2$ modulo $3$ we do not need to continue the weight reduction.\\

At a referee's suggestion, we include a picture describing the congruences involved in this step (see end of section 3.1 for notation):

$$
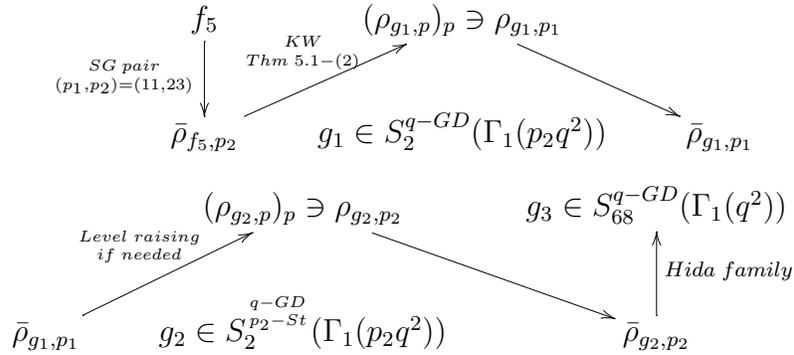
\xymatrix{
 f_5   \ar[d]_{SG \; pair \atop (p_1 , p_2)  = (11,23) } &   (\rho_{g_1,p})_p  \ni  \rho_{g_1 , p_1}  \ar[rd] &      \\
\bar{\rho}_{f_5, p_2}   \ar[ur]^{KW \atop Thm \; 5.1-(2)} &    g_1 \in  S_2^{q-GD}(\Gamma_1(p_2 q^2))     &     \bar{\rho}_{g_1 , p_1}    }
$$
$$
\xymatrix{ &   (\rho_{g_2, p})_p  \ni  \rho_{g_2 , p_2}  \ar[dr] &   g_3 \in S_{68}^{q-GD}(\Gamma_1(q^2))   \\
 \bar{\rho}_{g_1 , p_1}  \ar[ur]^{Level \;  raising \atop if \; needed}  &     g_2 \in  S_2^{q-GD \atop p_2-St}(\Gamma_1( p_2 q^2))      &     \bar{\rho}_{g_2, p_2}   \ar[u]_{Hida \; family } }
$$

Remark: $g_1$ has nebentypus of conductor $p_2$ and order $p_1$. \\
Having $g_3$, an application of WRGC (the technique explained in section 3.2) reduces to a newform $g_4$ of weight $32$ or $50$, and then another application of WRGC reduces the case of weight $50$ to weight $26$ or $38$. In any case, we finish this step with a newform $g_5 \in S_{26,32,38}^{q-GD}(\Gamma_1(q^2))$.

\subsection{Step 4}
Let us begin by giving the definition of Micro Good-Dihedral prime (as the referee points out, this notion has been used in [Di12a] but a formal definition is not given there):\\
{\bf Definition:} Let $f$ be a Hecke eigenform of weight $k \geq 2$ and level divisible by an odd prime $w$.
We say that $w$ is a Micro Good-Dihedral (MGD for short) for $f$ if the local type of $f$ at $w$ is
supercuspidal and the ramification of the compatible system of Galois representations attached to $f$
 at $w$ is given by a character (of the unramified quadratic extension of $\Q_w$) of order $s^\alpha \geq 7$
  where $s$ is a prime dividing $w+1$.\\

Remark: If the prime $s$ is odd, since $s$ can not divide $(w-1)
\cdot w$ it is clear that the nebentypus of $f$ (equivalently, the
determinant of the compatible system attached to $f$) can not ramify
at $w$. \\

 Thanks to the tricks in the previous step, we can
incorporate $43$ as a MGD prime, exactly as we did in [Di12a] to
incorporate the MGD prime $7$. Recall (cf. [Di12a]) that given a
prime $p$ that is the large one in a Sophie Germain pair, if we call
$\hat{p}_0 = (p-1)/2$, assumed prime, and we take another prime
$p_1$ dividing $p +1$, starting with an even weight $2 < k < p$, we
had devised a procedure to introduce $p$ as a MGD prime in the
level, in other words, to link through a chain of congruences the
given modular representation with another one, this time having $p$
as a supercuspidal prime in the level, with ramification given
by a character of order $p_1$, and having weight $2$. \\
This is explained in detail in  [Di12a], it exploits different known
cases of congruences between principal series, Steinberg, and
supercuspidal modular forms. First via a modular weight $2$ lift in
characteristic $p$ we introduce nebentypus at $p$ of order $\hat{p}_0$ (as
we did in the previous step with the pair $11$ and $23$), then in
characteristic $\hat{p}_0$ we transform this into Steinberg ramification
at $p$ (again, as in the previous step), and finally, moving to
characteristic $p_1$ we can create a non-minimal lift that transform
the type at $p$ into supercuspidal with ramification of order $p_1$.
As with Good-Dihedral primes, recall that this is useful because for
a compatible system of Galois representations attached to a modular form containing such a local parameter, the residual
representations in characteristics $\ell \neq p, p_1$ are going to
be irreducible (because they are so locally at $p$).\\
This is almost exactly what we do to introduce the MGD prime $p
=43$, using the auxiliary primes $p_0 = 7$ and $p_1 = 11$, to end up
with supercuspidal ramification of order $11$ at $43$. The only
difference is that $43$ is not Sophie Germain, i.e., $(43-1)/2 = 21$
is not a prime. But this is precisely why we have made Step 3. The
Sophie Germain trick, devised in the first place to introduce
Steinberg ramification at $p$, is based on the fact that, after
making a modular weight $2$ lift in characteristic $p$ the
nebentypus $\omega^{k-2}$ has order $\hat{p}_0$, a prime if we are in a
true Sophie Germain situation (observe that we have, as will always
be the case in this paper, even weights all through the weight
reduction process). Since we have forced the weight $k$ to be
congruent to $2$ modulo $3$ (this is what Step 3 is about), we can
work ``as if $43$ was Sophie Germain", the only important think is
that the nebentypus $\omega^{k-2}$ must have prime order, and this
is satisfied because since $k-2$ is divisible by $6$, this
nebentypus has order exactly $42/6 = 7$ (and, luckily, $7$ {\bf is} a
prime!), thus all the process of introducing the MGD prime can be
done
exactly as in [Di12a] thanks to the extra information on the weight.\\
Another key fact that allows us to do this step is that at the end of Step 3 the three possible values of the weight are all smaller than $43$ and larger than $2$: recall that the first thing we do in this step is to move to characteristic $43$, consider the residual representation and apply [KW], Theorem 5.1 - (2) to take a weight $2$ lift of this residual representation. The condition $ k < 43 $ is required to be in a Fontaine-Laffaille situation, so that the Serre's weight of the residual representation is known to be exactly this weight $2 < k < 43 $.  \\

We stress that all auxiliary primes must be taken to be greater than
$5$: this is the reason why we have not been able to find a
reasonable small true Sophie Germain prime that is good for our purposes.\\
As in previous sections, the Good-Dihedral prime $q$ with $q > t > B
> 68 > 43$ is enough to guarantee that all residual images in the above process are  $6$-extra large.\\
After these manipulations, we end up with a newform $f_6$ having
weight $2$ and level $43^2 \cdot q^2$, Good-Dihedral at $q$, and
having $43$ as a MGD prime. \\

At a referee's suggestion, we include a picture describing the congruences involved in this step (see end of section 3.1 for notation):

$$
\xymatrix{
g_5  \in  S_{26,32,38}^{q-GD}(\Gamma_1(q^2))   \ar[d] &   (\rho_{g_6,p})_p  \ni  \rho_{g_6 , 7}  \ar[rd] &      \\
\bar{\rho}_{g_5, 43}   \ar[ur]^{KW \atop Thm \; 5.1-(2)} &    g_6 \in  S_2^{q-GD}(\Gamma_1(43q^2))     &     \bar{\rho}_{g_6 , 7}    }
$$
$$
\xymatrix{
 &   (\rho_{g_7,p})_p  \ni  \rho_{g_7 , 11}  \ar[dr] &   f_6 \in S_2^{q-GD  \atop 43-MGD}(\Gamma_1(43^2 q^2))   \\
 \bar{\rho}_{g_6 , 7}  \ar[ur] ^{(*)}  &     g_7 \in  S_2^{q-GD \atop 43-St}(\Gamma_1(43q^2))      &     \bar{\rho}_{g_7 , 11}   \ar[u]_{KW \atop  Thm \; 5.1-(4)}  }
$$

We explain for the reader's convenience the arrow labelled $(*)$ (observe also that a similar situation already occurred in Step 3): since the nebentypus of $g_6$ is a character ramified at $43$ of order $7$, the residual representation $\bar{\rho}_{g_6,7}$ will be either unramified or with semistable (i.e., unipotent) ramification at $43$, and in both cases it was proved in [Di12a] that a lift that is Steinberg at $43$ exists. This will either be a case of Ribet's level raising (if the residual representations is unramified at $43$) or a minimal lift.

\subsection{Step 5}
At this step we are going to remove the Good-Dihedral prime $q$ from the level. \\
At the referee's request, let us explain why the prime $q$ is still a Good-Dihedral prime for $f_6$ as it was for $f_2$, i.e., with ramification given by a character of the same order $t$ and with respect to the same bound $B$. As in [KW] or in [Di12a] this follows from the definition of Good-Dihedral prime as long as one observes that in all congruences going from $f_2$ to $f_6$ the following conditions have been satisfied: all congruences are in characteristics smaller than $B$; whenever a prime $z$ has been added to the level in a congruence, it always was a prime satisfying $z < B$ and, finally,  in all congruences the local parameter at $q$ was not changed (all the lifts considered were minimal locally at $q$). Because of these conditions, we know that all congruences in the chain going from $f_2$ to $f_6$ had large residual image, this being forced by their local behavior at the Good-Dihedral prime $q$.\\
 In order to remove the Good-Dihedral prime from the level, we are first going to move to characteristic $t$, and then
 to characteristic $q$. Precisely because we are leaving the safety zone of
  characteristics up to $B$ we must be careful because in both cases, and
  in all the steps of the chain that follows, we need to check that the
  residual images are $6$-extra large: we don't have any longer the
  Good-Dihedral prime to ensure this.\\
Moreover, when reducing modulo $t$ we are losing the supercuspidal ramification
 at $q$, so we need other arguments even to justify that this residual representation,
 and those in the next steps of the chain, are irreducible.\\
On the other hand, we have also introduced the MGD prime $43$ in the level, and as
we already mentioned, as long as we avoid working in characteristics $43$ and $11$,
 the local information at $43$ is enough to guarantee that residual images are
 (absolutely) irreducible.
As the reader can easily check, in this step, and also in Steps 6,7 and 8, we
 are going to avoid these two characteristics, so residual irreducibility is guaranteed
 by the MGD prime in these four steps.\\
Also, ramification at $43$ has order $11$, and it gives an element of order $11$ also
 in the image of the projective representation, thus the exceptional cases in
  Dickson's classification of maximal subgroups of $\PSL(2, \F)$, where $\F$
  is a finite field of characteristic $p$, are ruled out for the images of
  the residual representations appearing in the chain, as long as the MGD
  prime stays in the level (again, this will be the case in Steps 5,6,7 and 8).\\
Therefore, in this step and the next three steps,  because of Dickson's classification
 and thanks to the MGD prime, in order to ensure that a residual image is large all
 that we will have to check is that it is not dihedral. \\
Moreover, since we are NOT going to work in characteristic $11$ during these four
steps, and using again the fact that we have an element of order $11$ in the
projective image, we know that if the residual image is large, it will also
 be $6$-extra large (this is a non-empty statement only for $p=13$).\\
This being said, let us proceed to remove the Good-Dihedral prime from the
level, checking that residual images are not dihedral.\\
We start with the newform $f_6$, we move to characteristic $t$ and we reduce
modulo $t$. \\
It is important to stress that, as happened in [KW], this is the exact moment in which the Good-Dihedral prime $q$ is lost: since the ramification at $q$ was given by a character of order $t$ and there are no elements of order $t$ in the multiplicative group of any finite extension of $\F_t$, it is clear that this character is killed when reducing mod $t$. Moreover, it is known (see [KW] and [Di12a] where the same situation occurred) that the residual mod $t$ representation will be either unramified or having semistable ramification at $q$ (and the two cases can occur, as proved in [KW], Thm 5.1 - (4)). \\
Suppose that the residual representation is dihedral, and let $K$
be the quadratic number field such that the restriction to $G_K$ becomes reducible.
 The only ramified primes for the residual representation are $t, q$ and $43$ (at most).
 Having reduced modulo $t$, the ramification at $q$ will be either trivial or semistable (unipotent).
 In either case, clearly $K$ can not ramify at $q$ (moreover, from the assumption of dihedral
 image we can conclude that the residual representation is unramified at $q$). On the other hand,
 the order of ramification at $43$ is odd, thus clearly $K$ can not ramify at $43$. We conclude
 that $K$ is the quadratic number field ramified only at $t$, in other words, that the residual
 representation is bad-dihedral.\\
But since the newform $f_6$ has weight $2$ and $t$ is not in its level, the Serre's weight of
 the residual representation is $k=2$, and then the residual representation can not be bad-dihedral
  because $t >3$ (cf. [Ri97]). The result of Ribet that we are applying here, a result that is part
  of the proof of ``generically large images for modular non-CM Galois representations"
   (cf.  [Ri85]), and is explained in detail in [Ri97] for the case of weight $2$, but holds
    in general with the same proof, is going to be required several times in the following sections,
     so let us record it here:\\

\begin{lemma}
\label{ribet} Let $\bar{\rho}$ be an irreducible representation in characteristic $p$,
 for a prime $p >2$, that is attached to some newform $f$. Suppose that the Serre's weight
  $k$ of $\bar{\rho}$ satisfies $2 \leq k \leq p+1$.
Suppose that $\bar{\rho}$ is bad-dihedral, or, more generally, that it is dihedral and
  induced from a quadratic number field that ramifies at $p$.\\
Then
it must hold: $p = 2k-1$ or $p =2k-3$. Moreover, in the first case the representation
is reducible locally at $p$, while in the second case it is, locally at $p$, induced from
a ramified character of the unramified quadratic extension $\Q_{p^2}$ of $\Q_p$.
\end{lemma}

This result of Ribet, that has been used many times in the literature (for example in [DM], or in [Kh]),
 can be easily proved by observing that in the case under consideration, the projective image of
  the inertia group at $p$ must have order $2$. The result follows easily from this and the
   definition of Serre's weight. The assumption on the size of $k$ is not serious, because
    it is known that by suitable twisting it is always satisfied (cf. [E]), and we always do apply
     this twisting to reduce to this situation, in all steps of the proof.\\

We have shown that the mod $t$ reduction of $f_6$ has $6$-extra large image. On the other
 hand, we know that it is either unramified or semistable locally at $q$. Moreover, in
 the first case, due to the existence of the lift given by $f_6$ which is supercuspidal
 at $q$ it is easy to see that the well-known necessary condition for Steinberg level-raising
  is satisfied. Thus, in both cases, doing level-raising if necessary, we consider a lift of
   this mod $t$ representation corresponding to a newform $f_7$ of weight $2$, level
    $43^2 \cdot q$, Steinberg at $q$.\\

Now we move to characteristic $q$, and we reduce modulo $q$. Using once again the MGD
prime $43$ in the level, we know that the residual image will be $6$-extra large provided
that it is not dihedral. Again, if it where dihedral, it can not be induced from a quadratic
 field $K$ ramified at $43$ because the order of ramification at this prime is odd, thus
  it should be bad-dihedral. But since $q$ is a Steinberg prime in the level, the Serre's
   weight $k$ of this mod $q$ representation is either $q+1$ or $2$, and this in any case
    by applying Lemma \ref{ribet} we know contradicts the fact that the representation
    is bad-dihedral. We thus conclude that the residual representation can not be
     dihedral, and thus that it has $6$-extra large image.\\

Since $f_7$ has weight $2$ and is Steinberg at $q$, we can consider the ordinary Hida family containing the $q$-adic
 representation attached to $f_7$. We now specialize this family at weight $q+1$ to get a newform $f_8$ which is ordinary at $q$ and congruent to $f_7$ modulo $q$,
 this newform $f_8$ of weight $q+1$ and level $43^2$ (see Section 3.3 for an explanation of why in a situation like this whenever one specializes to a weight congruent to $2$ modulo $q-1$ and larger than $2$ the level is prime to $q$).

\subsection{Step 6}
At this step we need again to perform the WRGC as in Step 2, in order to reduce
the (very large) weight $k = q+1$ of $f_8$ to a weight $k \leq 14$. \\
In this process, we are going to avoid characteristics $43$ and $11$, so that
 to ensure that the MGD prime $43$ is preserved. But we need to ensure that images
  are $6$-extra large in the whole process. In Step 2 this was ensured by the
   Good-Dihedral prime in the usual way (cf. [KW], [Di12b], [Di12a]), but now
    we only have the MGD prime in the level.\\
As explained in the previous section, thanks to the MGD prime, all that we need
 to check is that residual images are not dihedral. Moreover, in this process
 of WRGC, we will always work with mod $p$ Galois representations ramified only at $p$
 and at $43$, and since the order of ramification at $43$ is odd, the only possible
  dihedral case is the bad-dihedral case.\\
Thus, we have to perform WRGC, perhaps with some modifications with respect to the
 classical version given in Step 2, in such a way that we can guarantee that at each
  congruence the residual image is not bad-dihedral. \\
The process of weight reduction will follow a similar pattern as in Step 2:
we have a newform $f$ of weight $k$ and level $43^2$, such that $43$ is a MGD prime,
 and assuming that $k \geq  16$, we take $p$ to be the first prime larger than $k$, except
 for $k = 42$ where, since $43$ is forbidden, we are going to specify later what prime to
  take, and also for $k=32$ where in Step 2 we choose $p=43$ and here, again, we are going
   to have to choose another prime. \\
We reduce modulo $p$, and since it is easy to check (by hand for small values of $k$, and
because of the strong version of Bertrand's postulate used in Step 2 for $k > 36$) that we always have
 $p \neq 2k-1, 2k-3$ we apply Lemma \ref{ribet} and conclude that this mod $p$ representation
 is not bad-dihedral.\\
The problem is after doing the Galois conjugation: we can try formally to proceed exactly as
 in Step 2. With the same notation and same reasonings, we can take a modular weight $2$
  lift of the mod $p$ Galois representation, we know that the order $m$ of the nebentypus
   $\omega^{k-2}$ is greater than $6$, and we can take $t$ as before, the first number
    greater than $m/2$ relatively prime to $m$, and use it to Galois conjugate. The
    problem is that after doing so, the new Serre's weight of the mod $p$ reduction
    will be (after suitable twisting) equal to $k_1 = d t + 2$ or to $k_2 = p+1 - dt$
    (both smaller than the given $k$, as shown is Step 2), and in some cases it WILL
    BE THE CASE that some of these two values, call it $k'$, will satisfy $ p = 2 k' - 1$
    or $p = 2 k' - 3$, thus we can not apply Ribet's Lemma to rule out the bad-dihedral case. \\
We are going to solve this problem using the combination of two completely different tricks:
first, we will prove another Lemma that will allow us to eliminate the case $p = 2 k' - 3$,
 we will prove that, with our ramification conditions, such a bad-dihedral case can not happen.
 Then, once we are reduced to control the bad-dihedral case with $p = 2k' -1$, since this case
  can sometimes happen, we will isolate those cases where this equality is satisfied, and for
  such cases we will prove that there is another Galois conjugation, given by $t' > t$ the
   first prime to $m$ integer after $t$, such that by using this exponent to conjugate the
   new weights $k_1$ and $k_2$ are again smaller than the given $k$, and the equality $p = 2 k' -1$
   can not be satisfied by any of these two values (recall that we are assuming that it was
    satisfied by some of the corresponding values for $t$). In other words: it is not possible
     that using both conjugations, the one corresponding to $t$ and the one corresponding to $t'$,
      both give bad-dihedral cases. Thus, this extra conjugation will allows us to ensure that there
       is always a way to conjugate in order to reduce the weight, avoiding the bad-dihedral case. \\
The main difficulty is to show that $t'$ is still sufficiently small with respect to $m$, so as to
 imply that after this conjugation the new Serre's weights are smaller than the given one. \\
So, the first thing to do is to prove the result that allows us to get
 rid of the $p = 2 k' -3$ case. The following result is
  proved in [Kh] for the level $1$ case, but we are
  going to apply it in a more general situation,
  where it can be seen that it still holds applying the exact same reasoning:

\begin{lemma}
\label{winter} Let $\bar{\rho}$ be an irreducible representation in characteristic $p$,
for a prime $p >2$, that is attached to some newform $f$.  Suppose that the Serre's weight $k$
 of $\bar{\rho}$ satisfies $2 \leq k \leq p+1$. Suppose that for every prime $w$ different
  from $p$ where $\bar{\rho}$ is ramified, the order of the image of the inertia group at
  $w$ is odd. Then, the bad-dihedral case with $p = 2k-3$
  %, which corresponds to the local
  % at $p$ representation being induced from a ramified character of $\Q_{p^2}$,
   can not occur. \\
\end{lemma}
Remark: With the stronger assumption that the residual
representation ramifies only at $p$, this is Lemma 6.2 (i) in [Kh], and the
 proof given there extends to our case. For the reader convenience, let us
  briefly recall the argument. Suppose that the image of the representation
  is dihedral, i.e., we consider the projectivization $\mathbb{P}( \bar{\rho}  )$,
  and if we call $L$ the fixed field of its kernel, $L$ is a Galois number field
   with dihedral Galois group $D_{2t}$ of order $2t$, $t >1$. It follows from
   the hypothesis that the only quadratic number field $K$ that can be
   contained in $L$ is the quadratic number field unramified outside
   $p$ (in any other case, the image of inertia at some other
   prime would be even). From this, it follows that $t$ is odd,
    because for even $t$ dihedral groups are known to have more
     than one index two subgroup (even if in general only one of
      them is cyclic), therefore if $t$ were even $L$ would
      contain more than one quadratic field.\\
This is enough to rule out the case $p = 2k -3$, because in
 such a case locally at $p$ the image of the projective representation
  is a dihedral group of order $4$, but $D_{2t}$ can not contain
  a subgroup of order $4$ since we have shown that $t$ is odd.\\

Now we resume the process of WRGC: recall that we start with a representation of some weight $k >16$
 and level $43^2$, and we pick the first prime $p$ greater than $k$, except for $k = 32$
 and $k = 42$ where the value of $p$ will be specified later. In particular, we are always
  taking $p \neq 43$.\\
As in Step 2, the reduction process will be guaranteed by certain inequalities that hold
 for $k>36$ due to a strong version of Bertrand's postulate, and for $16 < k \leq 36$ and $k =42$
 it will be checked by hand (see Step 2). \\
Let us start with the case $ k > 36$, $k \neq 42$. We move to characteristic $p$ and reduce
 modulo $p$. As we already explained, Lemma \ref{ribet} implies that this residual
  representation is not bad-dihedral because $p$ is near $k$, more precisely: $p/k < 1.144$. \\
We take a modular weight $2$ lift (as usual, using [KW] Theorem 5.1 - (2) combined with a suitable M.L.T. such as Thm 4.1 in loc. cit.), whose nebentypus $\omega^{k-2} = (\omega^d)^a$ is a character
 of order $m = (p-1)/d > 6$, where $d = (p-1, k-2)$. Observe that the exponent $a$ is prime to $p-1$.
  By choosing some $t$ relatively prime to $m$, $m/2 < t < m-1$, which exists because $m > 6$,
   we consider the Galois conjugate of this representation having nebentypus $(\omega^d)^t$.
    As in Step 2, from the fact that $k-2$ is near $p-1$ (and thus dividing by $d$, the integer $a$ is
     near $m$) it follows that if $t$ is taken sufficiently near $m/2$, thus sufficiently
      away from $m$, the largest possible Serre's weight of the residual conjugated
       representation, $k_1 = d t +2$, will be smaller than $k$. In Step 2 we took
        $t$ to be the smallest possible number satisfying the above conditions. Now
         we need to be more careful, because we want to ensure that the residual conjugated
          representation is not bad-dihedral. We can apply Lemma \ref{winter}, because the
          only prime other than $p$ where this residual representation ramifies is $43$,
          and ramification at $43$ is of order $11$, thus we conclude that the bad-dihedral
           case with Serre's weight satisfying $p = 2 k' -3$ can not happen. Therefore,
           we only have to deal with the case $p = 2 k' - 1$. Observe that in all this
           process the Serre's weights are going to be even, because the determinant is
            unramified outside $p$, and modular Galois representations are odd. Then,
             in particular, if $p \equiv 1 \pmod{4}$, since $2 k' - 1 \equiv 3 \pmod{4}$
              we are already safe and we know that the bad-dihedral case can not happen.
              Later when we check by hand the weight reduction process for $k < 38$ this
               will be very important: whenever we choose the prime $p \equiv 1 \pmod{4}$
               then bad-dihedral cases can not occur (in this weight reduction).\\
Unfortunately, it is not clear that for any gigantic $k$ one can find a prime $p \equiv 1 \pmod{4}$
larger than $k$ sufficiently near $k$, so we still have to deal with the problem of bad-dihedral
 cases. The good thing is that we can restrict to the case $p \equiv 3 \pmod{4}$. With this
  restriction in mind, since $d = (p-1, k-2)$ is always even, we see that $m = (p-1)/d$ is odd. \\
Let us first use the same Galois conjugation that we used in Step 2, where $t$ is the first integer
larger than $m/2$ coprime to $m$. Since in our current situation $m$ is odd, this gives $t = (m+1)/2$. \\
Remark: in what follows, since we are restricted to the case $m$ odd, $m \geq 7$, we have $ \frac{1}{2} \phi(m) \geq 3$.
 This implies that there exists an integer $t'$ prime to $m$ with:
$$   m/2  <  t  = (m+1)/2 < t'  <  m-1  .$$
We conjugate using $t$, we compute the two possible values of the Serre's weight $k'$ of the
residual conjugated representation: $k_1 = dt +2$, $k_2 =  p + 3 - k_1$, and let us check what
are exactly the cases where one of these two values satisfies $p = 2 k' - 1$. First, recall
 that $k_1 + k_2 = p + 3$ and $k_1 > k_2$ (see Section 3.2).
  Thus, we have $ 2 k_1 - 1 > 2 k_1 - 3 > p$, and this implies that $k_2 = p + 1 - dt$, the smallest
  of the two values of the Serre's weight, is the only one that can give a bad-dihedral case.\\
So, let us suppose that we have:
$$ p = 2 k_2 - 1 = 2 (p+1 - dt) - 1  .$$
Or equivalently:
$$ p +1 - 2 d t = 0 .$$
Using  $p = md +1$, $m = 2t -1$, and expanding, we see that this is equivalent to:
$$ ((2t-1)d + 1 ) + 1 - 2 d t = 0 .$$
and this is equivalent to:
$$ d = 2 .$$
Therefore, we see that the (possibly) bad-dihedral case occurs exactly when $p \equiv 3 \pmod{4}$,
 $ t = (m+1)/2$ and $d = (p-1, k-2) = 2$.\\
We stress that except in this specific case, the process of WRGC is done exactly as in Step 2
using the value of $t$ as defined there, and what follows is an alternative Galois conjugation
 to be applied specifically in this case.\\
So, suppose that $p \equiv 3 \pmod{4}$, $d = (p-1, k-2) = 2$. Then, $m = (p-1)/2$ is odd, and
as in Step 2 we know that $m > 6$. \\
As we remarked, the following integer exists: let  $t' $ be the smallest integer prime to $m$
 satisfying $ t= (m+1)/2 < t' < m-1$. Let us use this value for conjugation instead of $t$.
 First of all, it is clear that using $t'$ we are not going to be in a bad-dihedral case,
 precisely because we are assuming that using $t$ we are in such a bad case. In fact, if both
 things happen at the same time we would have that in both cases $k_2$ satisfies: $p = 2 k_2 - 1$,
  thus the value of $k_2$ would be the same for both conjugated residual representations. Replacing
   by the values of these weights we would obtain the equality (recall that we are assuming that $d=2$)
$$ p + 1 - 2 t  =  p + 1 - 2 t' $$
contradicting the fact that $ t \neq t'$. \\
% In this argument we are using the fact that
%for $t $ and for  $t' $ (which is larger than $t$) we know that the largest Serre's
%weight $k_1$ is too big to correspond to the bad-dihedral case.\\
The key important fact, which remains to be checked, is that induction works also when one uses this new way of conjugating, i.e., that by conjugating using $t'$, the new Serre's weights are
smaller than the given weight $k$, and it is of course enough to check this for the largest
 of the two values, $k_1 = 2 t' + 2$.\\
In order to see this, we begin by proving the following elementary Lemma:
\begin{lemma}
\label{highschool} Let $m$ be an odd integer satisfying $m \geq 7$.\\
Let $p'$ be the smallest odd prime not dividing $m$. Then,  it holds:
$$ p' < 0.6 \cdot m $$
\end{lemma}
Proof: we divide the proof in four cases:\\
i) $3 \nmid m$: in this case, $p' =3$, and since $m \geq 7$ the Lemma follows
because $ 3 < 0.6 \cdot 7 = 4.2 $.\\
ii) $3 \mid m$ but $5 \nmid m$: in this case $p' = 5$ and since $m \geq 9$ the Lemma
follows because $5 < 0.6 \cdot 9 = 5.4 $. \\
 iii) $15 \mid m$ but $7 \nmid m$: in this case $p' = 7$ and since $m \geq 15$ the Lemma
  follows because $7 < 0.6 \cdot 15 = 9$.\\
iv) $105 \mid m$: Let us denote by $p_i$ the (positive) prime numbers in $\mathbb{Z}$, in its usual ordering. in this case if we have $p' = p_{r+1}$ it holds: $r \geq 4$ and $m$ is
 divisible by the $r-1 $ first odd primes $p_2 , p_3, ..... p_r $. Then, because by Bertrand's
  postulate $p' < 2 \cdot p_r$ we obtain:
$$ p ' < 9 \cdot p_r = 0.6 \cdot 3 \cdot 5 \cdot p_r = 0.6 \cdot p_2 \cdot p_3 \cdot p_r \leq 0.6 \cdot p_2 \cdot p_3 \cdot ...... \cdot p_r \leq 0.6 \cdot m $$
and this concludes the proof.
\qed \\
We have defined $t'$ to be the smallest integer that is prime to $m$ and larger than $(m+1)/2$. Then,
it is easy to see that if we define $p'$ as in the Lemma above, it holds:
$$ t' = (m + p')/2 .$$
Therefore, applying Lemma \ref{highschool} we conclude that $t' < 0.8 \cdot m$.\\
From this, we have: $k_1 = 2 t' + 2 < 1.6 m +2$, and since $m = (p-1)/2$ this gives:
$$ k_1 < 0.8 \cdot (p-1) + 2 = 0.8 \cdot p + 1.2 .$$
On the other hand, we are assuming that $k > 36$ in order to use, as in Step 2, the validity of a
 strong version of Bertrand's postulate that gives: $p  < 1.144 \cdot k$ (see formula (*) in section 3.2, where now we are calling $p$ the prime that was called $p_{n+1}$ there, this gives $p < 1.144 \cdot p_n$ which is stronger that what we want since $p_n < k$). Combining this inequality
  with the previous one we obtain:
$$ k_1 < 0.8 \cdot 1.144 \cdot k + 1.2 = 0.9152 \cdot k + 1.2 .$$
Also, since $k > 36$, it is easy to see that $0.9152 \cdot k + 1.2 < k$ and then we conclude that
 $k_1  < k$ giving the proof of the induction, for $k >36$ and $k \neq 42$.\\
It remains to see that for the remaining values of $k$
either $t$ defined as in Step 2  or  $t'$ as defined above makes the induction work. We will
 list all cases, starting from the special values $k = 42$ and $k= 32$. \\
In all cases  $p \neq 43$, $p \neq 2k-1$, and also $k_1, k_2 < k$. Moreover, since we are
 using $t'$ instead of $t$ precisely in those cases where $p \equiv 3 \pmod{4}$ and $d=2$,
  the equality $p = 2 k' - 1$ for $k' = k_1$ or $k_2$ is never satisfied:\\
$$ k =42,   p= 47:  d= 2,  m= 23, t' = 13, k_1 = 28, k_2 = 22 $$
$$ k = 32, p = 47: d=2, m= 23, t' = 13, k_1 = 28, k_2 = 22 $$
$$ k = 36, p = 37: d=2, m = 18, t = 11, k_1 = 24, k_2 = 16 $$
$$ k = 34, p = 37: d=4, m= 9, t = 5, k_1 = 22, k_2 = 18 $$
$$ k = 30, p = 31: d=2, m = 15, t' = 11, k_1 = 24, k_2 = 10 $$
$$ k = 28, p= 29: d=2, m=14, t = 9, k_1= 20, k_2 = 12 $$
$$ k = 26, p= 29: d= 4, m = 7, t= 4, k_1 = 18, k_2 = 14 $$
$$ k = 24, p = 29: d=2, m = 14, t = 9, k_1 = 20, k_2 = 12 $$
$$ k =  22, p = 23: d=2, m= 11, t' = 7, k_1 = 16, k_2 = 10 $$
$$k= 20, p=23: d= 2, m= 11, t' = 7, k_1 = 16, k_2 = 10 $$
$$ k = 18, p = 19, d= 2, m = 9, t' = 7, k_1 = 16, k_2 = 6 $$
$$ k = 16, p=17: d=2, m=8, t=5, k_1 = 12, k_2 = 8$$

This concludes the WRGC process, and we end up with a newform $f_9$ of level $43^2$ and
 weight $2 < k \leq 14$, $k$ even.\\

\subsection{Step 7}
This step is meant to reduce to cases where $k < 17 $ and divisible by $4$. This last condition
will be required in Step 8.\\
We start with $f_9$ and the idea is, as in Step 3, to unify the value of the
 weight using the Sophie Germain trick and then, since the weight obtained is too big,
  to perform some extra weight reductions, in a controlled way. This time we will apply
   Khare's method of weight reduction (cf. [Kh]). \\
To preserve the MGD prime $43$ in the level, we are going to avoid characteristics $11$
 and $43$ at this Step.\\
Once again, as explained at the beginning of Step 5, in order to ensure that residual
images are $6$-extra large it is enough to check that they are not dihedral.\\
So, let us first choose  the pair of Sophie Germain primes $23$ and $47$. We consider
a mod $47$ Galois representation attached to $f_9$. Since ramification at $43$ has odd order,
 the only possible dihedral case is the bad-dihedral case. Since the residual Serre's
 weight satisfies $2< k \leq 14 $ and $p=47$ it follows from Lemma \ref{ribet} that the
  representation can not be bad-dihedral, thus it can not be dihedral. \\
We take a modular weight $2$ lift, with nebentypus $\omega^{k-2}$ of
conductor $47$ and order $23$, $2<k \leq 14$, we move
 to characteristic $23$, and we reduce mod $23$. Observe that this residual representation
 has Serre's weight $2$ and is ramified at most at $47$, $23$ and $43$.
 Since we are reducing mod $23$ a character of order $23$, the ramification at $47$
  can either be trivial or unipotent.
  In any case, if we suppose that the residual representation has dihedral image, it must be
   unramified at $47$. Again, since $43$ is a MGD prime, the only possible dihedral case is
    the bad-dihedral case, which is ruled out by Lemma \ref{ribet} because $k=2$.\\
We take a lift of weight $2$ of this residual representation, corresponding to a newform
 that is Steinberg locally at $47$ (if the residual representation is ramified at $47$ this is just a minimal lift, and
 in case the residual representation is unramified at
  $47$, we are applying Ribet's level raising: as in step 3, it is easy to see that the necessary
  condition for level raising is satisfied). We move to characteristic $47$ and we
  reduce mod $47$. This residual representation can only have dihedral image if it
   has bad-dihedral image (the only ramified prime other than $47$ being the MGD
   prime $43$), which is ruled out again by Lemma \ref{ribet} because the Serre's
    weight is $k=48$ or $k=2$. \\
In any case, we can take an ordinary lift of this residual
representation of weight $48$ and level $43^2$: the specialization
at $k=48$ of the Hida family containing the weight $2$, Steinberg at
$47$ newform (as in Section 3.3, we know a priori that at weight
$48$ we specialize to a form of tame level $43^2$, but since $48$ is
congruent to $2$ modulo $47 - 1 = 46$ and is greater than $2$, we
conclude that the form has exact level $43^2$). Thus, we have been
able to reduce to a case of weight $48$. Even if $48$ is divisible
by $4$, it is too big for our purposes, so let us continue by
linking this representation
 with some other representation of smaller weight.\\
At this point, we apply a particular case of Khare's weight reduction: we pick a prime
 $p$ larger than $k$, we reduce modulo $p$ and take a modular weight $2$ lift, and we
 also select an auxiliary prime $r$ dividing $p-1$ so that by reducing modulo $r$ the
 order of the nebentypus changes. This way, one can take a lift of this mod $r$
 representation with a different type at $p$, in our case we are taking a lift
 that is minimal at $p$, and going back to characteristic $p$ the value of the
 residual Serre's weight will be different than the initial one (see [Kh], section 9, for details). \\
We have a newform of weight $48$ (and level $43^2$), and we choose $p=53$. We reduce
mod $53$ (as usual, the residual image is not dihedral because of Lemma \ref{ribet})
 and take a modular weight $2$ lift (whose existence follows from [KW], Thm. 5.1 - (2) combined
 with a suitable M.L.T. such as
  Thm. 4.1 in loc. cit.) with nebentypus $\omega^{46}$ of conductor $53$ and order $26$.\\
We move to  $r=13$ and we reduce modulo $13$ to kill a part of the
ramification at $53$ given by the character $\omega^{46}$. This
residual representation has Serre's weight $2$, and it possibly
ramifies at $13$, $43$ and $53$. Since ramification at $43$ has odd
order, and applying Lemma \ref{ribet}, we conclude that if the image
is dihedral, it must be the case that the representation is induced
from the quadratic
 number field ramified only at $53$. \\
Luckily, the MGD prime $43$ is a square modulo $53$, then we can
argue as with Good-Dihedral primes (this is the argument used to
obtain large images with Good-Dihedral primes, as in [KW], Lemma
6.3)
 to get a contradiction. For the reader's convenience, let us recall the argument.
  Let $K = \Q(\sqrt{53})$. Suppose that the residual representation is
  induced from a character of $G_K$, in which case its restriction to $G_K$ is reducible. Since
  $43$ is split in $K$, the image of the decomposition group at $43$ is contained in the restriction
   to $G_K$ of the representation, and this gives a contradiction since the first one is irreducible
    because $43$ is a MGD prime and the latter is reducible.  \\
This proves that the mod $13$ representation can not have dihedral image.\\
Before reducing modulo $13$, recall that we had a nebentypus at $53$
of order $26$. Thus, if we take a minimal lift of this mod $13$
representation (applying [KW], Thm. 5.1 - (3) combined with a
suitable M.L.T. such as Thm. 4.1. in loc. cit.), it will correspond
to a weight $2$ newform of level $43^2 \cdot 53$
but this time ramification at $53$ will be given by a character of order $26/13= 2$, namely, $\omega^{26}$.\\
We go back to characteristic $53$, and the possible values for the
Serre's weight of this residual representation are $26+2 = 28$, or
$53+3-28 = 28$, so $k=28$. This implies that the image is not
bad-dihedral even if $53 = 2 \cdot 28 -3$ , because we can apply
Lemma \ref{winter} (the only prime in the level is the MGD prime
$43$). Therefore,
 the residual image is is not dihedral.\\
So far, we have reduced to the case of a newform of weight $k=28$ and level $43^2$. Now we apply
 again Khare's weight reduction, taking $p=29$. We reduce mod $29$ (thanks to Lemma \ref{ribet}
  the residual image is not dihedral) and take a modular weight $2$ lift (using [KW], Thm. 5.1 - (2) combined with Thm. 4.1 in loc. cit.),
   whose nebentypus is
  $\omega^{26}$,
   a character of conductor $29$ and order $14$.\\
We move to $r=7$ and reduce modulo $7$ to kill a part of the
ramification at $29$ that was given by $\omega^{26}$, and obtain a
mod $7$ representation $ \bar{\rho} $ of weight $2$ ramified at $7$,
 $43$ and $29$. If we suppose that this representation is dihedral, induced from a character of
 $G_K$ for some quadratic number field $K$, then
  due to lemma \ref{ribet} we know that $K$ can not ramify at $7$, and since $43$ is a MGD prime we
  conclude
 that it must be $K = \Q(\sqrt{29})$.\\
It will take some effort to show that this is not possible: in fact
after reduction modulo $7$ ramification
 at $29$ is given by a character of order $14/7=2$, so a priori $ \bar{\rho}  $ could be induced from
 $K$. Notice, however, that the restriction $ \bar{\rho} |_{G_K} $
 will be unramified at $29$. \\
 Suppose that the image of $\bar{\rho}$ is dihedral. Let $\mu_1$, $\mu_2$ be the two mod $7$ characters whose sum is
isomorphic to the restriction to $G_K$ of $ \bar{\rho}  $. Each of
these two characters is possibly ramified at $43$ and $7$. We denote
by $\mathcal{O}_K$ the ring of integers of $K$ and by $v_1 , v_2$
the two places of $K$ above $7$. Let $\chi_{v_i}$, $i=1,2$, be the
composition $\mathcal{O}_K / v_i \xrightarrow{\sim} (\Z / 7
\Z)^\times \xrightarrow{\chi_7} \bar{\F}_7^\times $ where $\chi_7$
denotes the mod $7$ cyclotomic character. Choose a generator $\tau$
in $\Gal (K/\Q)$. Then $\chi_{v_1}^\tau = \chi_{v_2}$.\\
Since the image of $\mu_i$ is in characteristic $7$ and the
conductor of $\bar{\rho}$ at $43$ is $43^2$, by class field theory
we know that $\mu_i$ factors through $(\mathcal{O}_K / v_1 \cdot v_2
\cdot 43)^\times = (\Z / 7 \Z)^\times  \times (\Z / 7 \Z)^\times
\times
\F_{43^2}^\times $. \\
Let us call $\psi_i$ the character of conductor $v_i$ appearing on
the second component of $\mu_i$, so that $\mu_i$ can be written as:
$$ \mu_1  =  \chi_{v_1}^{i_1} \psi_1^{j_1} \phi_1, \;  \mu_2  =
\chi_{v_2}^{i_2} \psi_2^{j_2} \phi_2  ,$$ where $0 \leq i_1 , i_2,
j_1, j_2 \leq 6$ and $\phi_i$ is a character with conductor dividing
$43 \mathcal{O}_K$. \\
We fix an isomorphism of decomposition groups
$D_{v_1}\xrightarrow{\sim} D_7$, sending $\Frob_{v_1}$ to $\Frob_7$.
Note that $D_{v_2} = D_{v_1^\tau}  = \tau^{-1} D_{v_1} \tau$. In
this situation, by local information at $7$, since $\bar{\rho}$ has
Serre's weight $2$ we have: $ (\mu_1 \oplus \mu_2)|_{I_{v_1}} =
\chi_{v_1} \oplus 1$ and $(\mu_1 \oplus \mu_2)|_{I_{v_2}} =
 1 \oplus \chi_{v_2}$. Thus we must have
  $$\mu_1  = \chi_{v_1} \phi_1, \; \mu_2  = \chi_{v_2} \phi_2. $$
  Observe also that since the representation $\bar{\rho}$ was
  obtained from a newform with $43$ as a MGD prime in its level, its
  determinant (as well as the nebentypus of the newform) is unramified at $43$ (see remark after the definition of MGD prime in section
  3.4). Thus, since $K$ has class number
  $1$, we deduce that $\phi_2  \cdot \phi_1  = 1$. Thus, if $c$ denotes complex conjugation,
  we have $\phi_1 (c) = \phi_2 (c)^{-1} = \phi_2(c)$
   (because the image of $c$ has order at most $2$).\\
   On the other hand, $ \chi_{v_2}(c) = \chi_{v_1}^\tau(c) =
   \chi_{v_1}(c)= -1$. From this we conclude that $\mu_1 (c) =
   \mu_2( c)$ and therefore (recall that $K$ is real, so $c \in
   G_K$) that the determinant of $\bar{\rho}(c)$ is $1$,
   contradicting the fact that modular Galois representations are
   odd. \\
   This concludes the proof that the image of the representation
   $\bar{\rho}$ can not be dihedral. \\
   Remark: In a previous version of this paper, we gave a different proof which uses also some computations of ray class fields performed with Pari GP. In particular, these computations implied that the characters $\psi_i$ appearing in the first formula for the $\mu_i$ above have order at most $2$. The proof we have included above, proposed by the referee, does not use this extra information.\\

We now proceed as in the previous reduction: we take a modular
minimal lift of this residual representation (applying [KW], Thm.
5.1 - (3) combined with Thm. 4.1 of loc. cit.),
 observing that now the ramification at $29$ is given by the quadratic character $\omega^{14}$. \\
We go back to characteristic $29$, where the residual Serre's weight is $14 + 2 =16$ or $29 + 3 - 16 = 16$,
 and the residual image is not bad-dihedral (thus not dihedral) due to Lemma \ref{winter}. \\
We end up this step with a newform $f_{10}$ of weight $16$ and level $43^2$.\\

\subsection{Step 8}
%facil: introducir nebentypus en el $17$

This is a very elementary step, we have isolated this move just
because
of its conceptual importance. \\
Since the weight of $f_{10}$ is $16$, we can move to characteristic
$17$ and produce (applying [KW], Thm. 5.1 - (2) combined with Thm.
4.1 in loc. cit.) a congruence with a weight $2$ newform $f_{11}$ of
level $43^2 \cdot 17$, whose ramification at $17$ is given by the
character $\omega^{14}$ of order $8$, a character that is also the
nebentypus of this newform.\\
In this congruence, the residual representation has Serre's weight
$16$ and $p=17$, then an application of Lemma \ref{ribet} together
with the fact that $43$ is a MGD prime (where ramification has an
odd order) proves that the residual image is not dihedral. Due to
the MGD prime in the level, this is enough to conclude that the
residual image is $6$-extra large (see discussion at beginning of
Step 5).\\
This Step concludes with $f_{11}$, a weight $2$ newform with
nebentypus of order $8$ at $17$ and the MGD prime $43$ in the level.

\subsection{Step 9}
%magic lemma: modulo 11 se quita el MGD prime y la imagen es grande!
At this step we are going to remove the MGD prime $43$ from the
level of $f_{11}$, more precisely we will transform it into a
Steinberg prime
via a modulo $11$ congruence. \\
The problem is that since we are losing the MGD prime it is not
clear a priori that the residual image will be irreducible in this
congruence. Moreover, the space of newforms of weight $2$ and level
$43^2 \cdot 17$ seems to be too big for computations, so we would
not be able to check by hand that a mod $11$ representation attached
to $f_{11}$
is irreducible.\\
Luckily, the local information of the $11$-adic Galois
representation at the ramified primes $43$ and $17$ is enough to
deduce that the residual representation is irreducible, as proved in
the following:

\begin{lemma}
\label{magic} Consider the weight $2$ newform $f_{11}$, whose level
is $43^2 \cdot 17$, with nebentypus $\psi$ of order $8$ ramified at
$17$ and such that $43$ is a MGD prime with ramification at $43$
having order $11$. Then, any modulo $11$ residual representation
attached to $f_{11}$ is irreducible.
\end{lemma}

Proof: Let $\bar{\rho}$ be a residual representation in
characteristic $11$ attached to $f_{11}$. Suppose that $\bar{\rho}$
is reducible, and let us call $\mu_1$, $\mu_2$ the two characters in
the diagonal. Since $43$ was a supercuspidal prime, and ramification
at $43$ was given by a character of order $11$, it is well-known
that this residual representation will either be unramified or have
unipotent ramification at $43$. In any case, it is clear that the
characters $\mu_1$ and $\mu_2$ are unramified at $43$. From now on,
we consider the semi-simplification of the residual representation,
which is isomorphic to the direct sum $\mu_1 \oplus \mu_2$. Using
the information on the ramification at $17$, and the fact that the
newform is of weight $2$ and level prime to $11$, we conclude that
there are just two possibilities for this direct sum:
it must be either $\chi \psi \oplus 1$ or $\chi \oplus \psi$. \\
Here we are applying the usual principle that in a
Fontaine-Laffaille situation (residual characteristic $p=11 > k =
2$, $p$ not in the level) the residual Serre's weight equals the
weight of the modular form, combined with the fact that we know the
local inertial type
 at $17$ to be given by $\psi$. \\
 On the other hand, the $11$-adic lift of
$\bar{\rho}$ provided by $f_{11}$ is, locally at $43$, induced from
a character of the unramified quadratic extension of $\Q_{43}$, and
therefore the residual representation must satisfy the trace $0$
condition which is necessary for the existence of such a lift: in
terms of the characters $\mu_i$ the condition is:
$$ \mu_1 (43) + \mu_2 (43) = 0 $$
(recall that this is an equality in characteristic $11$).\\
Here we have used the fact that $\bar{\rho} (\Frob_{43})$ must
correspond to the action of $\Gal(\Q_{43^2}/\Q_{43})$, hence has the
form $\left(
        \begin{array}{cc}
          0 & 1 \\
          1 & 0 \\
        \end{array}
      \right) $
       giving trace $0$. \\

Observe that $ \chi(43) = -1 $.   We plug into the previous formula
the two possible values for the characters $\mu_i$ and in both cases
we obtain (again in characteristic $11$):
$$ \psi( 43 ) = 1 \quad \qquad (@)$$
Since $\psi$ is a character ramified at $17$ of order $8$, if we
take a prime $w$ which is a primitive root modulo $17$ the value
$\psi(w)$ will be an element of order $8$ in some extension of
$\F_{11}$. An easy computation shows that the order of $43$ modulo
$17$ is $8$, thus $43 \equiv w^2 \pmod{17}$ for some primitive root
$w$. Therefore,  $\psi(43) = \psi(w^2) = \psi(w)^2$ gives an element
of order $4$ in some extension of $\F_{11}$. But this clearly
contradicts (@), and this concludes the proof that $\bar{\rho}$ is
irreducible.\\
\qed

We have thus seen that the mod $11$ representation $\bar{\rho}$
attached to $f_{11}$ is irreducible. Let us check, using once again
Dickson's classification, that its image is large. First of all, the
projective image can not be an exceptional group because the image
of ramification at $17$ gives an element having, even after
projectivization, order $8$.\\
Since the Serre's weight is $2$, it follows from Lemma \ref{ribet}
that, if the image is dihedral, it must be induced from a quadratic
field $K$ ramifying at most at $17$ and $43$. But ramification at
$43$ of the residual representation is either trivial or unipotent
(thus in both cases of odd order), so $K$ can only ramify at $17$.
\\
To show that this case can not occur, we argue as in the proof of
Lemma \ref{ribet}. We are now considering ramification at a prime
different from the residual characteristic, but nevertheless the
same reasoning on dihedral groups applies: If we assume that the
image of $\bar{\rho}$ is dihedral induced from $K = \Q(\sqrt{17})$,
 using the fact that at the ramified prime $r =17$ the image of the
inertia group is given by a cyclic group of projective order greater
than $4$, we obtain a contradiction. \\
Thus, we conclude that the residual image can not be dihedral.
Therefore, from Dickson's classification we see that the image is
large. Moreover, since $\PGL(2, \F_{11})$ does not contain elements
of order $8$, we see (using a generator of the inertia group at
$17$) that the image is $6$-extra large.\\
As we have already remarked, the residual representation
$\bar{\rho}$ will have either trivial or unipotent ramification at
$43$. Moreover, in the unramified case,  the existence of the lift
given by $f_{11}$, which is supercuspidal at $43$, is enough to
ensure that Ribet's level raising can be applied to construct a lift
that is Steinberg at $43$. In fact, since $43 \equiv -1 \pmod{11}$,
the condition to apply Ribet's result reads
$$\bar{\rho}(\Frob_{43})
\equiv 0 \equiv 43 + 1 \pmod{11}, $$ and we have already stressed
during the proof of Lemma \ref{magic} that this is satisfied. In any
case, {\it soit} by Steinberg level-raising, {\it soit} by taking a
minimal lift, we conclude that there is a modular lift of
$\bar{\rho}$ corresponding to a weight $2$ newform $f_{12}$ of level
$43 \cdot 17$, which is
Steinberg at $43$ and has a nebentypus of order $8$ at $17$.\\

\subsection{Step 10}
What we do now is to move to characteristic $43$ and consider the
residual representation attached to $f_{12}$ in this characteristic.
In order to show that this residual representation is irreducible,
we perform a computation using MAGMA in the full space of newforms
of weight $2$ and level $43 \cdot 17$, with nebentypus given by any
character $\psi$ of order $8$ ramified at $17$ and trivial at $43$.
The output of this computation is that any residual representation
in characteristic
$43$ attached to any of these newforms is irreducible. \\
It took us only a few minutes to check this, because we used some
theoretical information to speed the computation. In fact, all we
did was to compute the Hecke polynomial $P_2(x)$ of the Hecke
operator $T_2$ in the above mentioned space (this is a polynomial with coefficients in $\Z$). In fact, and this will be clear in the explanation that follows, we only need
 to know the value of this polynomial modulo $43$, and this is what we computed. For the reader's convenience, we include some data of this polynomial. It has degree $256$ and:
 $$ P_2 (x )  \equiv     x^{256} + 8 x^{255} +  32 x^{254} + ..... + 33 x + 21  \pmod{43}  $$
 Then, we studied what
are the possible characters that would appear in the reducible case,
i.e., if the residual representation is assumed to be reducible, we
call $\mu_1$ and $\mu_2$ the characters in the diagonal, and using
local information we conclude that their sum will be either $\chi
\psi \oplus 1$ or $\chi \oplus \psi$, as in the previous step. This
follows from the fact that since $f_{12}$ has weight $2$ and is
Steinberg at $43$, its attached $43$-adic Galois representation is
ordinary and has residual Serre's weight $2$ or $43+1 = 44$.  We
evaluate this sum of characters at the prime $2$ in the two possible
cases and we consider $Q(x)$ to be the minimal polynomial of the
value obtained. Using the fact that $2$ is a square modulo $17$, we
see that the value $\psi(2)$ has order $4$. Thus, the two different
values for $Q(x)$ that we obtain are:
$$ (x-2)^4 - 1 , \quad \mbox{and}  \quad  (x-1)^4  - 2^4 $$
 To conclude,
we simply check using resultants that the polynomials $P_2(x)$ and
$Q(x)$ (we consider here the two possible values of the latter) do
not have any common root
modulo $43$ (observe that knowing only the value of $P_2(x)$ modulo $43$ is enough). \\
Knowing that the modulo $43$ representation is irreducible, let us
check that its image is $6$-extra large. Again, the ramification at
$17$ makes easy to see that the image can not be exceptional, and
for the dihedral case we just observe that since the Serre's weight
is either $2$ or $44$ it can not be dihedral corresponding to a
field $K$ ramifying at $43$ due to Lemma \ref{ribet}, and as in the
previous step we also eliminate the case where $K$ ramifies at $17$
because the image of the ramification group at $17$ contains an
element of projective order $8$.\\
We conclude from Dickson's classification that the image is
$6$-extra large.\\
We consider a modular lift of this residual representation,
corresponding to a newform $g$ of weight $44$ and level $17$, with
nebentypus at $17$ given by a character $\psi$ of order $8$. Observe
that both $f_{12}$ and $g$ are ordinary at $43$, they live in the
same Hida family.

\subsection{Step 11}
This step consists just on a computation performed using MAGMA on
the space of newforms of weight $44$ and level $17$, with any
nebentypus of order $8$ at $17$. Recall that $g$, the final output
of our chain, lives in this space. In agreement with some
expectations based on generalizations of Maeda's conjecture (work in
progress of P. Tsaknias jointly with the author), we conjectured
that this space
contains a unique orbit of Galois conjugated newforms.\\
This was confirmed by our MAGMA computation. This was computed by P.
Tsaknias using the MAGMA command NewformDecomposition, it took 14
hours in a CPU Intel Xeon E7-4850, 2GHz, with 192GB RAM.
 The computer used was
the server of the Number theory Group of The Department of
Mathematics of the University of Luxembourg. \\
For details on the computations involved and for the statement of
the generalization of Maeda's conjecture, see [DT].
\newline
Therefore, we have managed to link any given level $1$ cuspform $f$
with the orbit of $g$, and the chain constructed is safe, in the
sense that current A.L.T. (including the result in Appendix B) can
be applied in both directions, at any link, to the $5$-th symmetric
powers of the Galois representations considered. In particular, by
transitivity, together with the fact that Galois conjugation is a
valid move (so that we can move freely in the orbit of $g$), any
given pair of level $1$ cuspforms can be linked to each other in a
safe way. Thus, using our chain, we conclude that automorphy of
$\Symm^5(f)$ will hold for any level $1$ cuspform $f$  provided that
we can prove it for a single example.

\section{The base case for the automorphy of $\Symm^5(\GL(2))$}
For the base case, we start by considering an example of a weight
$1$ cuspform with projective image $A_5$ studied in [KiW]. This
example consists on a cuspform $f'_0$ of prime level $2083$, weight
$1$, and quadratic nebentypus $\epsilon = \omega^{1041}$ corresponding to $\Q(\sqrt{-2083}) / \Q$. We can consider this complex
representation as taking values on a finite field of characteristic
$2083$, and the obtained residual representation will have Serre's
level $1$ and Serre's weight $1042$ since $ \bar{\rho}_{f'_0 , 2083} |_{I_{2083}}  =
\left(
        \begin{array}{cc}
          \epsilon & 0 \\
          0 & 1 \\
        \end{array}
    \right) $, so we know that there is a
level $1$ cuspform $f_0$ of weight $1042$ such that this residual
representation is attached to it. Our goal is to prove automorphy
for $\Symm^5(f_0)$, by exploiting its congruence with
$\Symm^5(f'_0)$.\\
Using results of Kim and Wang (cf. [Kim], [W]), it is known that
$\Symm^5(f'_0)$ is automorphic and cuspidal. If we call $\rho$ the
complex Galois representation attached to $f'_0$, this uses the
identity: $$\Symm^5(\rho) \cong \Symm^2(\rho') \otimes \rho \qquad
\quad (*)$$ where $\rho'$ denotes the Galois conjugate
representation of $\rho$, which has trace defined over
$\Q(\sqrt{5})$. \\
But this is not enough to conclude residual automorphy of the
$2083$-adic Galois representation attached to $f_0$, in the sense
required for the application of A.L.T., because the automorphic form
attached to $\Symm^5(f'_0)$ is clearly not regular.\\
In order to solve this problem we first consider $\bar{\rho}$ and
$\bar{\rho'}$, the mod $2083$ representations obtained from $\rho$
and $\rho'$, and we replace them in formula $(*)$. Observe that
since $2083$ does not divide the order of the image the formula in
characteristic $2083$ is still an equality between irreducible
Galois representations. We now consider lifts of the residual
representations $\bar{\rho}$ and $\bar{\rho'}$ attached to cuspforms
of weights greater than $1$. More precisely, we know that we can
take a $2083$-adic Galois representation attached to $f_0$, a
cuspform of weight $1042$, as a lift of $\bar{\rho}$, and we can
take some $2083$-adic Galois representation attached to a weight $2$
newform $f_1$ minimally lifting $\bar{\rho'}$ (it exists because any
mod $p$ odd, irreducible, representation has a weight $2$ minimal
modular lift). Since $\bar{\rho'}$ also has ramification at $2083$
given by a quadratic character, it is clear that $f_1$ has level
$2083$ and that its
$2083$-adic Galois representations are potentially Barsotti-Tate (this modular lift is given for example by [KW], Thm. 5.1 - (2) combined with Thm. 4.1 in loc. cit.). \\
By plugging the representations attached to $f_0$ and $f_1$ into the
right hand side of formula $(*)$, and the one attached to $f_0$ on the
left hand side, we conclude that there is a modulo $2083$ congruence
between the Galois representations attached to $\Symm^5(f_0)$ and to
 $ \Symm^2(f_1) \otimes f_0 $. \\
 As in the work of Kim and Wang, we can deduce from known cases of
 Langlands functoriality (the $\Symm^2(\GL(2))$ case due to Gelbart and
 Jacquet in [GJ] and the $\GL(2) \times \GL(3)$ case due to Kim and Shahidi in [KS])
 that the latter tensor product is cuspidal automorphic (cuspidality follows from the criterion
  given in section 2 of [W]). Moreover, from our
 choice of the weights of $f_0$ and $f_1$ it is clear that its
 attached Galois representations are regular, and using the fact
 that for $f_0$ the prime $2083$ is not in the level and is in the
 Fontaine-Laffaille range and that $f_1$ is potentially
 Barsotti-Tate at $2083$ we see that the $2083$-adic Galois
 representation attached to the automorphic form $ \Symm^2(f_1) \otimes f_0
 $ is potentially diagonalizable (PO-DI). Here we are using the fact that being PO-DI is a property that is preserved by tensor products (this follows easily from remark 5 in page 530 of  [BLGGT], together with the fact that being PO-DI is a property that is preserved by arbitrary base change, as remarked in page 531 of loc. cit.) and also by symmetric powers (this, as we already explained at the beginning of section 3, follows from the same argument that gives preservation of being PO-DI by base change). \\
 On the other hand, we also know that the $2083$-adic Galois
 representation attached to $\Symm^5(f_0)$ is potentially
 diagonalizable (again, because $2083$ is in the Fontaine-Laffaille
 range for $f_0$). Finally, observe that the residual projective image in the
 congruence between $6$-dimensional Galois representations that we are considering is irreducible and
 isomorphic to $A_5$, thus clearly
 the residual representation
  will stay irreducible over any cyclotomic extension, thus it
  follows from the main result in  [GHTT]  (cf. Thm. 9 in loc. cit.) that its image is adequate
  even after restriction to a cyclotomic field. \\
  We have all the ingredients to apply the main A.L.T. in [BLGGT] (Thm 4.2.1 in loc. cit.) to
  deduce automorphy of $\Symm^5(f_0)$ from this mod $2083$
  congruence, because both $6$-dimensional Galois representations
  are regular and potentially diagonalizable, and the residual image
  satisfies the required condition. Therefore, the automorphy of $ \Symm^2(f_1) \otimes f_0
 $ implies that of $\Symm^5(f_0)$. This concludes the proof of the
 base case and, due to the results of the previous section, we also
 conclude automorphy of $\Symm^5(f)$ for any given level $1$ cuspform $f$.

\section{Application to base change}

A ``straightforward" concatenation of the safe chain constructed in
Section 3 with the process of ``killing ramification" (preceded by
ramification swapping) as used in section 3 of the paper [Di12a]
gives a new proof of base change for $\GL(2)$, for any newform of
odd level, this time without local conditions on the totally real
number field $F$. We do not even assume that $F$ is a Galois number
field. Thus, this proof works in much more general situations than
the one
given in [Di12a].\\
Moreover,  applying  the $2$-adic Modularity Lifting Theorem in [K-2] together with an adaptation
 of ideas from [KW] we can extend the proof of base change to newforms of arbitrary level by
 killing ramification at $2$ without increasing the weight.\\
We begin by explaining the proof for newforms of odd level, we will explain the case of even
 level at the end of this section. \\

 \subsection{The case of odd level}

 We start as in [Di12a] with a non-CM newform $f$ of odd level $N$ and weight $k \geq 2$
 and a given totally real number field $F$
  where we want to base change $f$. The assumption of being non-CM is harmless since it is well-known
   that base change holds for
   CM modular forms.

   \subsubsection{How to get rid of the level: reusing the method in ``fase
   uno" of [Di12a]}

   We go through all the process described in section 3 (called ``fase uno") of [Di12a], so as to safely link $f$ with
    a newform $f'$ of level $q^2$ and some weight $k'$, with $q$ as a Good-Dihedral prime with respect to a bound $B$ and $k' < B$.
     This prime $q$
    is added to the level in the same way as we did in section 3.1. \\
  Let us recall, for the reader's
 convenience, the main features of this chain.
 \\
 This safe chain (cf. [Di12a], section 3) consists of three main
 parts: \\
 (i) modulo an auxiliary prime $3 < r_0 < B$ chosen to be split in $F$ we move to a weight $2$ situation, and
 then working modulo a prime $t > B$
 a Good-Dihedral prime $q$ is added to the level
 (this step is similar to what we did in section 3.1, except that we do not conclude
 by increasing the weight again modulo $r_0$), this step ends with a newform $f_2$ of
 level $N \cdot r_0  \cdot q^2$ and weight $2$.\\
 (ii) ``ramification swapping" (cf. [Di12a], page
1024 to page 1027, line 6) is
 performed to safely link the newform $f_2$ to a newform $f_3$ of
 level $N' \cdot r_0 \cdot q^2$ and weight $2$, where  $N'$ is
 relatively prime to $N q$, in fact,  $N'$ is equal to the
 product of a set of odd primes $b_1 , b_2, ....b_z$ which are auxiliary
 primes chosen to be all split in $F$, larger than $N$, and smaller than the bound
 $B$ (the primes $b_i$, the prime $r_0$ and the bound $B$ are specified in advance, see [Di12a], page 1020). \\
 (iii) ``killing ramification" is performed: after certain manipulations that are required to ensure that the A.L.T.
  in [K-1] can be applied (in fact, such manipulations are no longer required: see
  the remark at the end of this subsection), by moving to each
 characteristic $b_i$ and $r_0$ and simply reducing modulo $b_i$ ($r_0$, respectively) and taking a
 minimal lift, all the primes $b_i$ and $r_0$ are removed from the level and
 we end with a newform $f'$ of level $q^2$ and some weight $k' < B$. \\
 This chain is safe in the sense that, after restricting to $G_F$,
at all the congruences a suitable A.L.T. applies in both directions.
Let us explain why this is true. \\
Concerning residual images, the
main point is that all residual images are large, as follows from
the introduction of a Good-Dihedral prime $q$ by the same arguments
we have used several times in section 3. It remains to prove that
being large is a property that is preserved when restricting to
$G_F$, because since large implies non-solvable, it is then
preserved when restricting to $G_{F(\zeta_p)}$, and also largeness
implies adequacy. But the proof that largeness is preserved when
restricting to $G_F$ follows exactly as in the proof of Lemma 3.2 in
[Di12a]: even if our current definition of large differs from the
one in loc. cit. and we want to include the case $p=3$ the {\bf
exact same argument}  applies (also note that now we are not
assuming that $F$ is Galois, but since the Galois closure $\hat{F}$
of $F$ is still totally real the proof applies to the restriction to
$G_{\hat{F}}$ thus a fortiori to the restriction to $G_F$). \\
Concerning the local conditions, for item (iii) recall that an
A.L.T. of Kisin was applied in [Di12a] precisely
 at this point,
 and it was checked that the conditions to apply this theorem over $F$ were satisfied (see also
 the remark at the end of this subsection
  for a simplified argument using an improved result). For item (i), the situation is
  as in Section 3.1, where the local conditions to apply the A.L.T. in [BLGGT] are satisfied over $\Q$, and
  also over $F$ (this is what we really need) because such conditions are preserved by arbitrary base change.\\
As for item (ii), the ramification swapping step, by the time the
paper [Di12a] was written, the results in this step were proved
under the assumption that $3$ and $5$ are split in $F$ (conditions
only required in case any of these primes is in the level). But this
is due to limitations inherent to the A.L.T. that were applied
there: if we use instead the A.L.T. from [BLGGT], and the one in our
Appendix B, we see that no local assumption is required. In fact, in
each congruence either both representations are ordinary or they are
both potentially diagonalizable, the latter case involving
potentially Barsotti-Tate and crystalline Fontaine-Laffaille
representations, as in previous sections (except when the residual
Serre's weight of a potentially Barsotti-Tate representation happens
to be $p+1$:  in this case a minimal crystalline lift of weight
$p+1$ is considered and this lift falls outside the
Fontaine-Laffaille range, but is known to be ordinary as proved in
[BLZ], thus again potentially diagonalizable). Thus, as it happened
in Section 3, at each congruence we are either in the ordinary case
or in the potentially diagonalizable case, and these conditions are
preserved by arbitrary base change,
 thus either Theorem 2.3.1 in [BLGGT] or the one in our Appendix B applies.\\
We conclude that thanks to the strength of current A.L.T., no local
assumption on the field $F$ (nor the fact that $F$ is Galois) is
required and still we can apply what is done in section 3 of [Di12a]
to safely link with a newform
 of level $q^2$, with $q$ as Good-Dihedral prime. \\
Remark: due to a recent improvement obtained in [HT] to the A.L.T.
of Kisin (the main result of [K-1]) that we have applied in the
killing ramification process, the safe chain constructed in section
3 of [Di12a] not only does not require any local
 assumption but also can be simplified. In fact, using the improved version
 (cf. [HT]) of this A.L.T., item (iii) above does
  not require any longer the tricks that were applied
  in [Di12a]: after completing item (ii) we can now go directly to each characteristic
  $b_i$ or $r_0$
   (observe that all these primes are greater than $3$),
   reduce modulo
   $b_i$ ($r_0$, respectively) and take a minimal lift, thus killing these primes from the level in a very elementary way.
   The only conditions required to apply this improved A.L.T. are
   that the residual image has to be irreducible even after restriction to $G_{F(\zeta_{b_i})}$ (since large implies
    non-solvable, we know that this holds), and the residual characteristic has to
   be larger than $3$ and split in $F$, which holds by construction.

  \subsubsection{Using the chain from section 3, starting at Step 2}

At this point, we have already reduced the proof of base change, at least for
 odd level, to the construction we performed in Section 3: once $f$ is safely
 linked to a newform $f'$ of level $q^2$, we can start from $f'$ and apply
 all the construction described in Section 3 except that we start now from
 Step 2 (Section 3.2) because the Good-Dihedral prime is already in the
 level. \\
 There is only one caveat: the newform $f'$ has level $q^2$ and
 weight $k' \geq 2$, and the problem is that we do not know how to
 eliminate the possibility that when finishing
  the killing ramification process we end up with $k' = 2$, and the condition $k' > 2$
 is required to go on  in Step 2 (see Section 3.2).
This can easily be remedied by acting as if the Good-Dihedral prime
$q$ in the level were not a Good-Dihedral prime, then choosing and
adding another prime $q'$ much bigger such that it is Good-Dihedral
for all characteristics up to some bound $B'$ larger than $q$. Then,
removing $q$ from the level in two steps (as done in Section 3.5),
we know that we will end up with a newform of weight $q+1 >2$ and
level $q'^2$, good-dihedral at $q'$ (compare with section 3.5, where
after removing $q$ from the level in two steps, first working mod
$t$ and then mod $q$, we end up with the newform $f_8$ of weight
$q+1$). It is clear that in these steps the A.L.T. in [BLGGT] can be
applied to check that automorphy propagates well among the
restrictions to $G_F$ (as in the previous subsection, the key points
are that the required local conditions and largeness of residual
images are both preserved by base changing to $F$). We conclude that
by a simple iteration of the argument we can assume ``without loss
of generality" that after the completion of the killing ramification
step the weight obtained is greater than $2$. \\
 Thanks to this trick, we can go through the process in Sections 3.2 to 3.11 and thus we obtain a safe chain
  linking $f'$, and therefore also $f$, to the single orbit space
   computed at section 3.11. Observe that all that we are taking
   from
    section 3 is also safe in this new context, because by base
    changing to $G_{F(\zeta_p)}$
     large residual images stay large (as we explained in the previous subsection) and are thus adequate,
     and the local conditions in the A.L.T.
    that we are applying (Theorem 2.3.1 in [BLGGT] and the one in our Appendix B) are preserved by base change. \\

    \subsubsection{A CM form as base case}

As in the proof of automorphy for $\Symm^5(\GL(2))$, we see that all
that is required to complete the argument is a base case, i.e.,
 if there is a newform $f_0$ satisfying base change, then since both $f$ and $f_0$ can be safely
 linked to the same single orbited space, then base change propagates from one to the other.
 A priori, we were only considering non-CM forms, and this was useful to have generically
 large residual image. However, we will pick as $f_0$ a CM form, which is very convenient
  since for such an $f_0$ arbitrary base change is well-known, let us see how we can make
  one single congruence to transform such an $f_0$ into a non-CM form $f'_0$ in such a way
   that this congruence is safe in our context so that existence of base change for $f_0$ implies the same for $f'_0$.
    After having done this, base change is
   propagated from $f'_0$ to $f$ by the argument already explained, and in particular in
    a way that all congruences involved will have large residual image (the chain
    linking $f$ to $f'_0$ is just a concatenation of the two chains that link both
     forms with the single orbit space described in section 3.11).
So let us pick a suitable base case. It is easy to see that given any form
$f_0$ attached to a CM elliptic curve $E$ of odd conductor, if we
call $K$ the field such that $E$ has CM by an order of $K$, $f_0$
can be base changed to any totally real number field $F$ (as it is
well-known, since $K$ can not be contained in $F$ the restrictions
to $G_F$ of the $\ell$-adic Galois representations attached to $f_0$
are absolutely irreducible). Also, given $F$ we easily see that
there is a prime $p > 5$ (depending on $F$) not in the level of
$f_0$ such that the residual mod $p$ representation attached to
$f_0$ (it is of course dihedral) will be irreducible even after
restriction to $F(\zeta_p)$ (we just need to ensure that
$\hat{F}(\zeta_p)$ does not contain $K$, where $\hat{F}$ denotes the
Galois closure of $F$, but it is easy to see that this can only fail
for finitely many primes).\\
Then, modulo such a prime $p$ the first thing that we do is to apply
level-raising to add a Steinberg odd prime $w$ to the level (keeping
the weight equal to $2$). At this mod $p$ congruence, it follows
from the main result in [GHTT] that the residual image restricted to
$F(\zeta_p)$ is adequate (and both representations are
Barsotti-Tate), thus it follows from the A.L.T. in our Appendix B
that modularity over $F$ propagates well. After introduction of the
Steinberg prime we already have a non-CM weight $2$ newform $f'_0$
of odd level with generically large residual images (cf. [Ri85]),
and we can proceed to introduce the Good-Dihedral prime to the level
and go on with all the known process (ramification swapping, killing
ramification and so on) to safely link it with the single orbit
space in Section 3.11, thus also with the given newform $f$. Except
at this first step where the residual image is dihedral (but not
bad-dihedral), all the residual representations
in the rest of the chain will be large. \\
This concludes the the proof of base change for odd level newforms
where, thanks to the fact that the local conditions in the A.L.T. in
[BLGGT] and in Appendix B are preserved by base change, we do not
need to impose any local
condition on the field $F$. \\

\subsection{The case of even level}

Let us now focus on the case of a newform $f$ of even level. As usual, we assume that
$f$ is a non-CM form. The idea to prove base change for such a newform is to eliminate
 the prime $2$ from the level as in [KW], relying on the $2$-adic M.L.T in [K-2], thus
  reducing the proof to the odd level case, already solved. \\
First of all, since we are going to need the $2$-adic M.L.T. over $F$, let us observe
that we are going to apply Theorem 0.9 in [K-2], but with the local condition at the
prime $2$, numbered as condition (3) in loc. cit., replaced by the condition:\\
(3') Let $v$ be any prime of $F$ dividing $2$ such that the representation $\rho$ restricted
to the decomposition group at $v$ is ordinary. Then the $2$-adic representation $\rho_h$,
which is attached to a Hilbert newform $h$ and is congruent to $\rho$ modulo $2$, is also
 ordinary locally at $v$.\\
The fact that Theorem 0.9 in loc. cit. holds with this modification can be seen by inspecting
 its proof (the same happens with the M.L.T. in odd residual characteristics proved in [K-BT],
 where a version of the main theorem with this more general condition is recorded), it could
 have been stated this way. \\
We proceed with $f$ as we did in the odd level case, the only difference is that at some point
 we are going to eliminate the prime $2$ from the level. The right moment for this is after
 the swapping ramification process (applied to all ODD primes in the level), and before starting
 the killing ramification. At this point the newform,  let us call it $f'$, contains in its
  level the prime $2$, plus a set of relatively large auxiliary primes $b_i$ and $r_0$, and a very
  large prime $q$, the good-dihedral prime, where ramification has order $t$, another
   very large prime number. Observe also that at this point of the argument $f'$ has
   weight $2$ (cf. [Di12a]).  \\
Because of the Good-Dihedral prime, we know that if we work modulo a small odd
characteristic the residual image will be adequate. Also, as in [KW] or [Di12a],
 the prime $q$ is supposed to satisfy certain conditions with respect to the prime $2$,
  and the prime $t$ is taken greater than $5$, and using this it follows (cf. [KW])
  that if we work modulo $2$ the residual projective images will be non-solvable and
  not isomorphic to $A_5$.\\
Let us also assume, as in [Di12a], that the prime $q$ is split in the Galois closure
 $\hat{F}$ of $F$, this way after restricting to $G_{\hat{F}}$  a $2$-adic Galois
 representation we know that the residual projective image still is absolutely
  irreducible and contains an element of order $t$, and from this we conclude,
  as in [Di12a], Lemma 3.4, that the image of this restriction is again
  non-solvable and not isomorphic to $A_5$. We easily conclude   that the
  same holds for the image of the restriction to $G_F$. \\
Let us now indicate the moves, taken from [KW], that will allow us to link
 the compatible system attached to $f'$ with another corresponding to a newform
  of odd level, in such a way that, after restricting to $G_F$, modularity propagates
   well from the latter to the former.\\
We divide in two cases:\\
(i) the $2$-adic Galois representation attached to $f'$ is potentially Barsotti-Tate, \\
(ii) the $2$-adic Galois representation attached to $f'$ is potentially semistable of weight $2$.\\
Let us start with case (ii). In this case, we can twist the representation
 by a suitable character in order to reduce to the semistable case,
  so let us assume that the $2$-adic representation attached to $f'$
  is semistable. As in [KW], section 9, we reduce modulo $3$ to change
   the local type at $2$, preserving the weight and the ramification
    at all other primes (this is based on Theorem 5.1 - (3) of loc. cit.):
     since $3$ is not in the level of $f'$ (the primes $b_i$ and $r_0$ are large primes)
      this is just a congruence between Galois representations that are both
       Barsotti-Tate at $3$, and because of the good-dihedral prime the residual image is large,
thus adequate, therefore we know from the A.L.T. in our Appendix B
that the new Galois representation we are creating is also modular,
attached to certain newform $f_2$. This new Galois representation
constructed in [KW] has a different type at $2$, $f'$ was Steinberg
while $f_2$ is principal series at $2$ (ramification being given by
characters of order $3$). It is important to observe that the
$2$-adic Galois representation attached to $f_2$ is potentially
Barsotti-Tate. Moreover, it is shown in [KW] that the residual
modulo $2$ Galois representation attached to $f_2$ has Serre's
weight $2$. Therefore, applying
 again Theorem 5.1 of loc. cit. (recall that we know that the projective residual
  image is non-solvable and not an $A_5$) we can take a lift of this mod $2$ representation
   corresponding to a newform $f_3$ of weight $2$ and odd level. More precisely,
   the lift is constructed using Theorem 5.1 - (1) in loc. cit., and its modularity
    follows from Theorem 4.1 in loc. cit. \\
Since $f_3$ has odd level, we know that it can be base changed to $F$. Let us now explain
why the chain we have just constructed linking $f$ to $f_3$ is safe, in the sense that it
allows to propagate modularity for the restrictions to $G_F$, backwards. It is enough to
concentrate on the part of the chain that goes from $f'$ to $f_3$, because we have already
 explained (when dealing with the odd level case) that  the part linking $f$ to $f'$ is safe
 (the fact that now the prime $2$ is in the level does not affect the argument). \\
The restrictions to $G_F$ of the modular Galois representations being considered have residual
 images that are adequate, even after restriction to $G_{F(\zeta_p)}$,  when the residual
  characteristic is $p=3$, and non-solvable and not projectively $A_5$ when $p=2$. For the
   mod $3$ congruence between $f'$ and $f_2$,
we easily see that the A.L.T. in our Appendix B applies over $F$.
Thus,
 it remains to check that in the mod $2$ congruence modularity over $F$ propagates from $f_3$
  to $f_2$, applying the modification of Theorem 0.9 in [K-2] discussed above.
  Conditions (1) (residual modularity and non-solvable residual image) and (2) (potentially
  Barsotti-Tate at primes above $2$) of this theorem are satisfied by the restriction to
   $G_F$ of the $2$-adic representation attached to $f_2$, and concerning condition (3'),
    if we assume that the $2$-adic Galois representation attached to $f_2$ satisfies this
     potential ordinarity condition, then it is known that this can only happen if it is
      nearly-ordinary at $2$ (cf. [H], Prop. 3.3), but this implies that the residual
       mod $2$ representation is ordinary, thus that the Barsotti-Tate representation
       attached to $f_3$ is ordinary. \\
We conclude that the $2$-adic M.L.T. of Kisin applies, thus that modularity propagates
well over $F$, from $f_3$ to $f$.\\
Let us now treat case (i). This time we move directly to
characteristic $2$ and we reduce modulo $2$. If the residual Serre's
weight is $2$, we take a lift corresponding to an odd level, weight
$2$ modular form $f_2$ (the existence of which follows from Theorems
4.1 and 5.1 - (1) of [KW]),
 and we see as in the previous paragraph that the $2$-adic M.L.T. of Kisin
  allows to propagate modularity over $F$ from $f_2$ to $f'$. If the residual
   Serre's weight is $4$, we take a lift, whose existence is guaranteed by
   Theorems 4.1 and 5.1 - (2) of loc. cit., corresponding to a newform $f_2$ of weight $2$
    and even level whose level is strictly divisible by $2$ and the prime $2$ is Steinberg.
     Since $f_2$ falls in case (ii), we have already shown that it can be base changed
     to $F$. It remains to check that the $2$-adic M.L.T. of Kisin can be applied to
     propagate modularity over $F$ from $f_2$ to $f'$, and again condition (3')
     is the only non-trivial one. But since $f_2$ is semistable at $2$ and has weight $2$, it is
      also ordinary at $2$, so condition (3') holds automatically. \\
We conclude that, in any case, we can construct a safe chain linking $f$ with
an odd level newform, and in particular, that $f$ can also be base changed to $F$.\\

%Moreover, by using results of Khare-Wintenberger (a $2-3$ trick that allows to reduce to a potentially Barsotti-Tate situation at $2$) and Kisin (a $2$-adic M.L.T. for potentially Barsotti-Tate representations) it might be %possible to kill the ramification at $p=2$ in order to extend the above proof of base change to the case of even level modular forms. We plan to explore this possibility in the near future.\\

\section{Base change for $\Symm^5(\GL(2))$}

For any newform $f$ of level $1$ we have shown automorphy of
$\Symm^5(f)$, let us see that by combining this with the base change
result in the previous section we can
also deduce that $\Symm^5(f)$ can be base-changed to any totally real number field $F$.\\
This follows by considering the same chain constructed in section 3, and the proof of automorphy
 for the base case given in section 4, and checking that all the construction can be
 base changed to $F$. The fact that the safe chain also works well over $F$ is
 automatic, if we argue as in the previous section: the local conditions for
 the A.L.T are preserved by base change, and after restriction to $G_F$ a
 residual image that is $6$-extra large remains $6$-extra large. \\
Concerning the base case, the projective residual image $A_5$ is not
changed by restriction to $G_F$ (again, you can use the fact that the projective
 image of the restriction to $G_{\hat{F}}$  is not trivial because the
  representation is odd and the Galois closure $\hat{F}$ of $F$ is totally
  real, therefore since $A_5$ is simple the projective image of the
   restriction to  $G_{\hat{F}}$ does not change, thus clearly for
    the restriction to $G_F$ the same holds). The rest of the proof
    of automorphy for the base case, this time over $F$, goes exactly
    as in section 4, just notice that $f_0$ and $f_1$ can be base changed
    to $F$ because of the result in the previous section, and that the
    known cases of Langlands functoriality applied to conclude that
    $\Symm^2(f_1) \otimes f_0 $ is automorphic are also known to hold
    when $f_0$ and $f_1$ are Hilbert newforms over $F$.

\section{Final Comments}

Around the same period of time where the results in this paper were obtained (summer 2012), Clozel and Thorne developed a different, yet related,
approach to prove functoriality of some symmetric powers of
$\GL(2)$, using an A.L.T. proved by Thorne that works in residually
reducible situations. With this technique they managed to deduce
automorphy of $\Symm^5(\GL(2))$ and $\Symm^7(\GL(2))$ for classical
modular forms of arbitrary level, and also for many Hilbert
newforms. Their approach shares several similarities with ours (for example, the use of Good-Dihedral primes, the use of the main
A.L.T. from [BLGGT] and the use of the
criteria taken from our Appendix A to
guarantee adequacy of certain residual images). The
 fundamental difference is that they have  new, powerful tools, at
 their disposal (an A.L.T. working in the residually reducible case, together with new ``level-raising"
 results)
 so they can consider a case of small characteristic that we were forced to avoid.
  The interested reader should consult [CT].\\

\section{Bibliography}

[BLGGT]   Barnet-Lamb, T., Gee, T., Geraghty, D., Taylor, R., {\it
Potential automorphy and change of weight}, Ann. Math. {\bf 179}
(2014)  501-609
\newline
[BLZ] Berger, L., Li, H., Zhu, J., {\it Construction of some families of $2$-dimensional crystalline representations}, Math. Ann. {\bf 329} (2004),  365-377
\newline
[CHT] Clozel, L., Harris, M., Taylor, R., {\it   Automorphy for some $\ell$-adic lifts of automorphic mod $\ell$ representations},
Pub. Math. IHES {\bf 108} (2008), 1-181
\newline
[CT] Clozel, L., Thorne, J., {\it Level-raising and symmetric power
functoriality, II}, Ann. Math., to appear
\newline
[Di09] Dieulefait, L., {\it  The level $1$ case of Serre's conjecture revisited}, Rendiconti Lincei - Mat. e Appl. {\bf 20} (2009),  339-346
\newline
[Di12a] Dieulefait, L., {\it Langlands base change for $\GL(2)$},
Ann.   Math. {\bf 176} (2012) 1015-1038
\newline
[Di12b] Dieulefait, L., {\it Remarks on Serre's modularity conjecture},
Manuscripta Math. {\bf 139} (2012) 71-89
\newline
[Di12c] Dieulefait, L., {Automorphy of $m$-fold tensor products of
$\GL(2)$}, preprint; available at www.arxiv.org
\newline
[DM] Dieulefait, L., Manoharmayum, J., {\it  Modularity of rigid Calabi-Yau threefolds over $\Q$}, Calabi-Yau Varieties and Mirror Symmetry. N. Yui, J. Lewis  et al. (eds.). Fields Institute Communications vol. 38;  (2003) AMS, 159-166
\newline
[DG] Dieulefait, L., Gee, T., {\it Automorphy lifting for small
$l$}, Appendix B to this paper
\newline
[DT] Dieulefait, L., Tsaknias, P., {\it On the generalized Maeda's conjecture and its relation to local types}, in preparation
\newline
 [Ge] Geraghty, D., {\it Modularity lifting theorems for ordinary
Galois representations}, preprint
\newline
[E] Edixhoven, B., {\it The weight in Serre's conjectures on modular forms}, Invent. Math. {\bf 109} (1992) 563-594
\newline
[GL] Gao, H., Liu, T., {\it A note on potential diagonalizability of
crystalline representations}, Math. Ann., to appear
\newline
[GK] Gee, T., Kisin, M., {\it  The Breuil-Mezard conjecture for
potentially Barsotti-Tate representations}, preprint
\newline
[GJ] Gelbart, S.,  Jacquet, H., {\it A relation between automorphic representations
of $\GL(2)$ and $\GL(3)$}, Ann. Scient. Ecole. Norm. Sup. {\bf 11} (1979) 471-542
\newline
[G] Guralnick, R., {\it  Adequacy of representations of finite
groups of Lie type},  Appendix A to this paper
\newline
[GHTT]  Guralnick, R., Herzig, F.,  Taylor, R.,  Thorne, J.,  {\it
Adequate subgroups}, appendix to ``On the automorphy of $\ell$-adic
Galois representations with small residual image" by J.Thorne, J.
Inst. Math. Jussieu {\bf 11} (2012) 855-920
\newline
[H] Hida, H., {\it A finiteness property of abelian varieties with potentially ordinary good reduction}, J. Amer. Math. Soc. {\bf 25} (2012) 813-826
\newline
[HT] Hu, Y.,  Tan, F., {\it The breuil-mezard conjecture for
non-scalar split residual representations}, arXiv:1309.1658
[math.NT], 2013
\newline
 [Kim] Kim, H., {\it An example of
non-normal quintic automorphic induction and modularity of symmetric
powers of cusp forms of icosahedral type}, Invent. Math. {\bf 156}
(2004), 495-502
\newline
[KS] Kim, H., Shahidi, F., {\it    Functorial products for $\GL(2) \times \GL(3)$ and
the symmetric cube for $\GL(2)$}, Annals of Math. {\bf 155} (2002) 837-893
\newline
[KiW] Kiming, I., Wang, X., {\it Examples of $2$-dimensional odd
Galois representations of $A_5$-type over $\Q$ satisfying the Artin
conjecture}, in ``On Artin's Conjecture for odd $2$-dimensional
representations", G. Frey (Ed.), LNM 1585, (1994) Springer-Verlag
\newline
[K-BT] Kisin, M., {\it Moduli of finite flat group schemes, and modularity}, Annals of Math. {\bf 170} (2009) 1085-1180
\newline
[K-1] Kisin, M., {\it The Fontaine-Mazur conjecture for $\GL_2$},
J.A.M.S {\bf 22} (2009) 641-690
\newline
[K-2] Kisin, M., {\it  Modularity of $2$-adic Barsotti-Tate representations}, Invent. Math. {\bf 178} (2009) 587-634
\newline
[Kh] Khare, C., {\it     Serre's modularity conjecture: The level one case}, Duke Math. J. {\bf 134} (2006)  557-589
\newline
[KW] Khare, C., Wintenberger, J-P., {\it Serre's modularity conjecture (I)}, Invent. Math.  {\bf 178} (2009) 485-504
\newline
[Ri85] Ribet, K., {\it On $\ell$-adic representations attached to modular forms. II}, Glasgow Math. J. {\bf 27} (1985) 185-194
\newline
%[Ri90] Ribet, K., {\it Raising the levels of modular representations}, S\'{e}minaire de Th\'{e}orie des Nombres, Paris 1987--88, 259-271, Progr. Math. {\bf 81}, Birkhauser (1990)
%\newline
[Ri97] Ribet, K., {\it Images of semistable Galois
representations}, Pacific J.  Math. {\bf 181} (1997)
%[RS] Ribet, K., Stein, W., {\it Lectures on Serre's conjectures}, Arithmetic Algebraic Geometry (Park City, UT, 1999), 143-232, IAS/Park City Math. Ser. {\bf 9}, A. M. S. (2001)
%[R3] Ribet, K., {\it On modular representations of ${\rm Gal}(\overline{\Q}/ \Q)$ arising from modular forms}, Invent. Math. {\bf 100} (1990) 431-476
\newline
[Sa] Savitt, D., {\it On a Conjecture of Conrad, Diamond, and Taylor}, Duke Math.
J. {\bf 128} (2005)141-197
\newline
[W] Wang, S. {\it On the symmetric powers of cusp forms on $\GL(2)$
of icosahedral type}, Int. Math. Res. Not. (2003), 2373-2390.

\end{document}